\documentclass[final,leqno]{siamltex}
\usepackage{amsmath}
\usepackage{amsfonts}
\usepackage{mathrsfs}
\usepackage{amssymb}
\usepackage{pgfplots}
\usepackage{tikz}
\usepackage{pgfplots}
\pgfplotsset{compat=1.5}
\usepackage{subfig}
\usepackage{graphicx}
\usepackage{bm}
\usepackage{longtable}
\newtheorem{example}{}[section]
\newtheorem{remark}{\bf Remark}[section]
\newcommand{\norm}[1]{\lVert #1 \rVert}
\newcommand{\tnorm}[1]{\ensuremath{\left| \! \left| \! \left|} #1 \ensuremath{\right| \! \right| \! \right|}}
\newcommand{\mc}[1]{\mathcal #1}
\newcommand{\e}[1]{\epsilon #1}
\newcommand{\T}{\mathcal T}
\newcommand{\J}{\jmath}
\title{CONVERGENCE OF AN ADAPTIVE MIXED FINITE ELEMENT METHOD FOR GENERAL SECOND ORDER LINEAR ELLIPTIC PROBLEMS}

\author{Asha K. Dond\thanks{Department of Mathematics, Indian Institute of Technology, Bombay
        ({\tt asha@math.iitb.ac.in}).}
       \and Neela Nataraj \thanks{Department of Mathematics, Indian Institute of Technology, Bombay
                ({\tt neela@math.iitb.ac.in})}
                \and Amiya Kumar Pani \thanks{Department of Mathematics, Indian Institute of Technology, Bombay
                                ({\tt akp@math.iitb.ac.in})}
                }

\begin{document}
\maketitle
\slugger{mms}{xxxx}{xx}{x}{x--x}

\begin{abstract}
The convergence  of an adaptive mixed finite element method for general second order linear elliptic problems defined on simply connected bounded polygonal domains is analyzed in this paper. The main difficulties in the analysis are posed by the non-symmetric and  indefinite form of the problem along with the lack of the orthogonality property in mixed finite element methods. The important tools in the  analysis are {\it a~posteriori} error estimators, quasi-orthogonality property and quasi-discrete reliability established using  representation formula for the lowest-order Raviart-Thomas solution  in terms of the Crouzeix-Raviart solution of the problem. An adaptive marking in each step for 
the local refinement is based on  the edge residual  and  volume residual terms of the {\it a~posteriori} estimator. 
Numerical experiments confirm the theoretical analysis.
\end{abstract}


\pagestyle{myheadings}
\thispagestyle{plain}
\markboth{Convergence of AMFEM}{}
\section{Introduction}
The general second-order linear elliptic PDE  
on a  simply connected bounded  polygonal Lipschitz domain  $\Omega \subset {\mathbb R}^2$ with boundary $\partial \Omega$ reads: For given right-hand side $f\in L^2(\Omega)$, seek $u$ such that
\begin{eqnarray}\label{eq1}
 \mathcal{L}u :=  -\nabla \cdot (\mathbf A \nabla u+u {\mathbf b})+ \gamma ~ u=f  \hspace{.5cm}\mbox {in}~\Omega, ~~~u=0~\hspace{.5cm}\mbox{on}~\partial\Omega.
\end{eqnarray}
The coefficients ${\bf A,b},\gamma$ are all piecewise smooth and the symmetric matrix ${\mathbf A}$ is positive definite and uniformly  bounded away from zero.\\
The flux variable ${\bf p}  =-({\mathbf A\nabla u+u {\mathbf b}})$ and ${\mathbf b}^*= \textbf A^{-1}{\bf b}$ 
 allow to  recast  (\ref{eq1}) as a  first-order system  
\begin{eqnarray}\label{eq3}
\begin{array}{lll} 
\mathbf A^{-1}{\bf p} +u {\mathbf b}^* + \nabla u =0  \; ~{\text{and}}~~
{\text{div}}~{\bf p}+\gamma \: u=f \;  {\rm in } \; \Omega. 
\end{array}
\end{eqnarray}
The mixed formulation  
seeks  $({\bf p},u)\in H(\text{div},\Omega) \times L^2(\Omega) $ such that 
\begin{eqnarray} \label{eq4}
\begin{array}{llll} 
(\textbf A^{-1}{ \bf p} + u {\bf b}^*,{\bf q})-({\rm div}~ {\bf q}, u)=0 \qquad \mbox {for all}\; {{\bf q}} \in H(\text{div},\Omega),\\
({\text{div}}~ {\bf p}, v )+(\gamma \: u,v)=(f,v) \qquad \qquad \mbox{for all}\; v \in L^2(\Omega).
\end{array}
\end{eqnarray}
Here and throughout the paper, $H(\text{div},\Omega)=\{{\bf q}\in L^2(\Omega; \mathbb R^2) : \: \text{div}~ {\bf q} \in L^2(\Omega) \} $ and 
$L^2(\Omega; \mathbb R^2)$ denotes the space of $\mathbb R^2$-valued  $L^2$  functions defined over the domain $\Omega$.
The existence and uniqueness of the mixed solution for elliptic problems have been proved in \cite{Brezzi,CC}. 

 The study of analysis of the adaptive finite element methods (AFEM) is an essential component  of the adaptive process. Various  {\it a~posteriori}  error estimators are  reviewed in \cite{oden}, and the references therein.
 The marking strategies, 
 convergence and optimality 
 are well established for  the adaptive conforming finite element methods in literature \cite{Car_apo1,quasi, dof, morin1, morin2,  veeser, steven}. 
For the Poisson problem, the convergence and optimality   
 have been established for the adaptive nonconforming FEM  \cite{bek_non,CH,car_non_cgt} 
and for the adaptive mixed FEM   \cite{bek_mixed,car_cgt2,car_cgt1,cascon,chen_cgs}.  
The recent article `Axioms of adaptivity' {\cite{axiom}} provides a general framework  to optimality of adaptive schemes.

The non-symmetric and indefinite second order elliptic equations  with conforming, nonconforming mixed FEM have been discussed in various articles \cite{Brenner-Scott, Brezzi,CC,Chen,chen_hoppe,daglas,schatz1,schatz}. These articles discuss the existence and uniqueness of the solution with {\it a~priori}  error estimates. 
 {\it A posteriori}  error estimates and its  convergence for conforming FEM for  general second order 
 linear elliptic PDEs have been achieved using contraction of the sum of energy error plus oscillation  
 in  \cite{nochetto} and  the quasi- optimality in \cite{prac}.  
{\it A posteriori}  error estimates  and quasi-optimal convergence of the adaptive nonconforming FEM   
have been obtained in  
  \cite{chen_hoppe}. 
To the best of our knowledge, we have not come across any work which discusses the 
convergence and optimality of the adaptive mixed finite element method (AMFEM) for non-symmetric and indefinite elliptic problems. 
The main challenges, the lack of orthogonality in MFEM and the non-symmetric form of equation  are addressed in this work.
Also as the flux variable  ${\bf p}$ involves  $u$ explicitly,  
 the analysis of variable $u$  becomes inevitable for the analysis of the flux ${\bf p}$. 
In this paper, the main contributions are summarized as:
\begin{itemize}
\item for the adaptive algorithm, the  marking strategy in each step for 
the local refinement is proposed  based on the comparison of the edge residual term and the volume residual terms of the {\it a~posteriori} estimator,
\item {\it a~posteriori} error estimator, quasi-orthogonality property and 
quasi-discrete reliability results are derived with the help of the representation formula for the lowest-order Raviart-Thomas solution  in terms of 
the Crouzeix-Raviart solution of the problem,
\item the contraction property is shown for the linear combination of 
the sum of errors in  ${\bf p}$ and $u$, 
the edge residual  estimator  and the volume residual estimators,
\item  the convergence and the quasi-optimality results are achieved, under the assumption of  small initial mesh-size $h_0$.
\end{itemize}
An outline of the  paper is as follows.
 Section 2 introduces  notations and the adaptive algorithm for the mixed finite element method.
Section 3 describes some auxiliary results necessary for the convergence analysis. 
The contraction property
 and the quasi-optimal convergence of the adaptive mixed finite element method  are established in Section 4.
The numerical experiments  are presented in Section 5. Appendix I summarizes the constants used in the article and their interdependencies.

Here are some notations used  throughout the paper. An inequality $A\lesssim B$ abbreviates $A\leq CB$, where $C>0$ is a mesh-size independent constant that depends only on the domain and 
the shape of finite elements;
 $A\approx B$ means $A \lesssim B \lesssim _{â€¢} A$. Standard notation applies to Lebesgue and Sobolev spaces and $\norm{\cdot}$ abbreviates
 $\norm{\cdot}_{L^2(\Omega)}$ with   $L^2$ scalar product $(\cdot,\cdot)$. For a vector ${\bf q}=(q_1,q_2)\in H({\rm div},~ \Omega)$, $\norm{{\bf q}}:=(\norm{q_1}^2+\norm{q_2}^2)^{1/2}$. Let $\norm{\cdot}_T $  and $\norm{\cdot}_E $ 
 denote  $\norm{\cdot}_{L^2( T)} $ and $\norm{\cdot}_{L^2(E)} $ respectively. $H^m(\Omega)$ denotes the  Sobolev  space of order $m$ with  norm given by $\norm{\cdot}_m.$  
\section{AMFEM algorithm}
This section discusses notations and  the  adaptive algorithm. \\
Let $\mathcal{T}_h$ be a regular triangulation  the domain 
${\Omega}\subset {\mathbb R}^2 $ into triangles such that 
$\cup_{T\in \mathcal{T}_h} T =\overline{\Omega}.$
Let $\mathcal{E}_h $ be the set of all edges in $\mathcal{T}_h$ and let ${\mathcal{E}_h}({\partial \Omega})$ be the 
set of all boundary edges in $\mathcal{T}_h$. 
Further, let $\text{mid}(E) $ denote the midpoint of the edge $E$ and $\text{mid}(T)$ denote the centroid of the triangle $T.$ 
 The set of edges of the element $T$ is denoted by $ {\mathcal{E}} (T)$ and let 
 $h_T:=$ diam$(T)$ for $T \in \mathcal{T}_h$. Define $h_{\mathcal T} \in P_0(\mathcal T_h)$, as a piecewise constant mesh-size function such that   
$h_{\mathcal T}|_T:= h_T$ for all $T \in \mathcal T_h$. Let $h:= \max_{T\in \mathcal {T}_h} h_T$ and  
$h_E$ be the length of the edge $E\in \mathcal{E}_h $. 
For any edge $E $, $\nu_E $ is  the unit normal vector exterior
to $T$ and $\tau_E$ is the unit tangential vector along $E$.
Let $\Pi_h$ be the $L^2$ projection onto  $P_0(\mathcal{T}_h)$ and define $osc_h(f):=\norm{h_{\mathcal T}(1-\Pi_h)f},$ where
$$
P_r(\mathcal{T}_h)=\{v \in L^2(\Omega): \forall T\in \mathcal{T}_h, v|_T\in P_r(T) \}.
$$
Here, $P_r(T)$ denotes the algebraic polynomials  of total degree at most $r \in \mathbb {N}$ as functions 
on the triangle $T\in \mathcal{T}_h.$
%
The jump of ${\bf q}$ across $E$ is denoted by $[{\bf q}]_E$; that is, for  two  neighboring triangles
$T_{+}$ and $T_{-},$
$$
[{\bf q}]_E(x) :=({\bf q}|_{T_{+}} (x)-{\bf q}|_{T_{-}}(x))~\text{  for } ~x \in E =\partial T_{+} \cap 
\partial T_{-}.
$$
The sign of $[{\bf q}]_E$ is defined using the convention that there is a fixed orientation of $\nu_E$ pointing 
outside of $T_{+}.$ 
The patch $\omega_E$ denote the union of elements that share a common edge $E$. 
The piecewise gradient $\nabla_{NC}: H^1(\mathcal {T}_h)\longrightarrow L^2(\Omega; {\mathbb R}^2)$ acts as
$ \nabla_{NC}  {v}|_T = \nabla v|_T\;\;\;\text{for all} ~T \in \mathcal{T}_h.$ 
The broken Sobolev norm $\tnorm{\cdot}_{NC}$ abbreviates $(\mathbf{A} \nabla_{NC} ~\cdot,\nabla_{NC} ~\cdot)^{1/2}$.\\
The non-conforming Crouzeix-Raviart (CR) finite element space  with respect to the triangulation $\mathcal T_h$ reads
\begin{eqnarray*}
&& CR^1(\mathcal T_h): =\{v \in P_1(\mathcal{T}_h): v~ \text{ is continuous in all midpoints mid($E$) of edges} ~E \in \mathcal E_h\},\\
&& CR^1_0(\mathcal T_h ):=\{v \in CR^1(\mathcal T_h):v(\text{mid} (E))
=0 ~~\text{for all} ~E \in \mathcal E_h ({\partial\Omega}) \}.
\end{eqnarray*}
The lowest-order Raviart-Thomas space with respect to the triangulation $\mathcal T_h$  reads 
\begin{align*}
RT_0(\mathcal{T}_h):=\{{\bf q}\in H(\text{div},\Omega): & \: \forall T \in \mathcal{T}_h 
~\exists {\bf c} \in \mathbb{R}^2~ \exists d \in \mathbb{R}~ \; \forall {\bf x} \in T, ~{\bf q}({\bf x})={\bf c}+ d ~ {\bf x}\\
&\text{ and} \; ~\forall E\in \mathcal{E}_h, [{\bf q}]_{E} \cdot \nu_E=0\} .  
\end{align*}
\subsection{ Algorithm}
 The standard structure of an adaptive algorithm is successive loops 
 $$\mbox{ SOLVE $\rightarrow$ ESTIMATE $\rightarrow$ MARK $\rightarrow$ REFINE}$$
 on different levels of the triangulation.\\
 In the step { SOLVE}, the discrete mixed finite element problem (RTFEM) for \eqref{eq4} defined by: seek $ ({\bf p_h}, u_h) \in RT_0(\mathcal {T}_h)\times P_0(\mathcal T_h)$ such that
\begin{eqnarray}
 (\mathbf{A}^{-1}{\bf p_h}+ u_h {\bf b}^*,{\bf q}_h)-(\text{div}~{\bf q}_h,u_h)=0 \qquad 
 \text{for all}\, {\bf q}_h\in RT_0(\mathcal {T}_h), \label{eqna1}\\
 (\text{div}~ {\bf p_h},v_h)+ (\gamma u_h, v_h)=(f_h,v_h) \qquad  \text {for all} \, v_h \in P_0(\mathcal{T}_h), \label{eqna2}
\end{eqnarray} 
is solved.  
Recall that, $f_h:=\Pi_hf $ is the $L^2$-projection of $f$ onto $P_0(\T_h)$. \\ 
The step {ESTIMATE} consists of computation of an {\it a~posteriori} error estimator. Here, {\it a~posteriori} error estimator is a combination of the edge estimator $\eta_h$ and
the volume  estimator $\mu_h$, that is,
\begin{equation}\label{esti_fin}
\norm{{\bf A^{-{1}/{2}}(p-p_h)}}^2+\norm{u-u_h}^2\leq C_{rel} (\eta_h^2+\mu_h^2),~~~~{\rm where}
\end{equation}
 \begin{align}
   \eta_h^2&:=\sum_{E\in \mathcal E_h}\eta_h^2(E)~ {\rm with} ~ \eta_h^2(E):=\norm{h_{E}^{1/2}[{\bf A^{-1}p_h+b^*}u_h]\cdot \tau_E}^2_E ;\label{eta}\\
  \mu_h^2&:=\sum_{T \in \mathcal T_h }\mu_h^2(T)~~{\rm with}~\nonumber\\
  &~~~\mu_h^2(T):=osc_h^2(f)_T
  +\norm{h ~{\rm div}~{\bf p_h}}^2_T+\norm{h ({\bf A^{-1}p_h }+u_h{\bf b^*})}^2_T. \label{mu}
\end{align}
 The {\it a~posteriori} estimate (\ref{esti_fin}) is  derived in Section 3, Theorem 3.2.\\
The step {MARK} consists of the two alternatives {\bf (A)} and {\bf (B)} which depend on the computable quantities $\eta_h$ and $\mu_h$ and a positive parameter
$\kappa$. \\
  {\bf Case (A)}: if $\mu^2_h\leq \kappa \eta_h^2 $, compute the minimal set of edges $\mathcal M_h \subset \mathcal E_h $ such that 
 \begin{equation}\label{A}
  \eta_h^2 \leq \eta_h^2(\mathcal M_h) ~~{\rm with}~~0<\theta_A<1.
   \end{equation}
 {\bf  Case (B)}: if $\mu^2_h> \kappa \eta_h^2 $,
compute the minimal set of triangles $\mathcal M_h \subset \mathcal T_h$ such that 
\begin{equation}\label{B}
\theta_B \mu_h^2 \leq \mu_h^2(\mathcal M_h)~~{\rm with}~~0<\theta_B<1. \end{equation}
Here $\theta_A$, $\theta_B$ and $\kappa$ are the parameters of the marking criteria and will be chosen  appropriately.\\
 Newest vertex bisection (NVB) algorithm  \cite{biven,steven} is  applied for refining the marked  edges  or elements and generate a new regular triangulation in  {REFINE} step.
Note that to maintain the conformity of the triangulation, some additional edges and elements may also need refinement.
\begin{remark} Instead of separate marking (\ref{A})-(\ref{B}) defined in MARK,  one could use a collective marking, that is, compute the minimal set of triangles $\mathcal M_h \subset \mathcal T_h$
\begin{equation}
\theta_A (\eta_h^2 + \mu_h^2 ) \leq (\eta_h^2+\mu_h^2)(\mathcal M_h) ~~{\rm with}~~0<\theta<1.
\end{equation}
\end{remark}
\section{Auxiliary results}
This section discusses some important results required  for the convergence analysis which are the {\it a~posteriori} error estimator, error and estimator reduction properties.\\
The  nonconforming finite element method (NCFEM) for (\ref{eq1}) seeks
 $ {u}_{CR}\in CR^1_0(\mathcal T_h)$ such that 
\begin{equation}\label{non1}
(\mathbf{A}\nabla_{NC}{u}_{CR}+ {u}_{CR}  {\bf b} ,\nabla_{NC} v_{CR})+(\gamma{u}_{CR},v_{CR})= (f, v_{CR}),{\forall} ~ v_{CR} \in CR^1_0(\mathcal T_h).
\end{equation}
{\bf Representation of RTFEM Solution via NCFEM} \cite{CC,Marini}: 
The coefficients ${\bf A,b},\gamma$ are all piecewise constants.
Now the auxiliary discrete problem is to seek  $ u_h^N\in CR^1_0(\mathcal T_h)$ such that
\begin{equation}\label{eqna3}
 (\mathbf{A}\nabla_{NC} u_h^N+{u}_h {\bf b},\nabla_{NC} v_{CR})+(\gamma ~{u}_h,v_{CR})=
 (f_h,v_{CR}),\text{for all}\, v_{CR}\in CR^1_0(\mathcal T_h),
 \end{equation}
where for $T \in \mathcal{T}_h$ 
\begin{eqnarray} 
&& \displaystyle  u_h({\bf x})=\left(1+\frac{S(T)}{4} \gamma \right)^{-1} \left(\Pi_0 u_h^N+\frac{S(T)}{4}f_h \right)  \;\; \text{for}~ {\bf x} \in T,\label{umt}\\
&& ~~S(T)= \displaystyle{\int_T} ({\bf x}-\text{mid}(T))\cdot \mathbf{A}^{-1}({\bf x}-\text{mid}(T)) \,d{\bf x}.\label{s}
\end{eqnarray}
Then, the solution ${\bf p}_h$ of the mixed finite element method formulation (\ref{eqna1})-(\ref{eqna2}) satisfies
\begin{equation} \label{pmt}
 {\bf p}_h({\bf x})=-\left(\mathbf{A}\nabla _{NC} u_h^N+ u_h{\bf b} \right)+\left(f_h-\gamma  u_h \right)
 \frac{\left({\bf x}-\text{\rm mid}(T)\right)}{2} \text{ for}~ {\bf x}\in T.
\end{equation}
The well-posedness of (\ref{non1}) and (\ref{eqna3}) and the equivalence of (\ref{eqna1})-(\ref{eqna2}) with (\ref{eqna3}) is discussed in \cite{CC}.
\begin{lemma}\label{lm8}
Let $u_h^N$ and $({\bf p_h}, u_h)$ solve (\ref{eqna3}) and (\ref{eqna1})-(\ref{eqna2}), respectively. Then, it holds
\begin{eqnarray}
&& \norm{\nabla_{NC} u_h^N}\leq \norm{\mathbf{A}^{-1}{\bf p_h}+u_h{\bf b}^*},\label{5.16}\\
&&\norm{{\rm div}~{\bf p_h}} = \norm{f_h-\gamma u_h} \lesssim \norm{f_h}+\norm{\mathbf{A}^{-1}{\bf p_h}+u_h{\bf b}^*}.\label{5.17}
\end{eqnarray}
\end{lemma}
{\em Proof}. 
From (\ref{pmt}), 
 $${\bf A}^{-1} {\bf p_h}+u_h{\bf b}^*=-\nabla_{NC} u_h^N+\left(f_h-\gamma u_h \right){\bf A}^{-1}
\displaystyle{ \frac{\left({\bf x}-\text{\rm mid}(T)\right)}{2}}.$$
Since 
$\left((f_h-\gamma  u_h)\left({\bf x}-\text{mid}(T)\right)/2, \nabla_{NC} u_h^N \right) =0,$
the Pythagoras theorem yields
\begin{equation*}
 \norm{{\bf A}^{-1} {\bf p_h}+u_h{\bf b}^*}^2=\norm{\nabla_{NC} u_h^N}^2
 +\norm{\left(f_h-\gamma u_h \right)\mathbf{A}^{-1} 
\displaystyle{ \frac{\left({\bf x}-\text{\rm mid}(T)\right)}{2}}}^2. \hspace{1cm}
\end{equation*}
Hence, (\ref{5.16}) holds. A use of triangle inequality with  (\ref{umt}) and (\ref{5.16})  implies (\ref{5.17}). \hfill${\Box}$
\\
The  following theorem  is on   {\it a~posteriori}  error estimates of  ${\bf e_p:= p -p_h} $~and $e_u:=u-u_h$ 
the proof of which is obtained by minor  modifications in the proof of Theorem 5.5 in \cite{CC}. However, for the sake of completeness, a short proof is given below. 
\begin{theorem}\label{5.5}({\it A~posteriori} error estimate)
 ~Let $u \in H^1_0(\Omega) $ be the unique weak solution of (\ref{eq1}) and 
  let $({\bf p_h},u_h) $ be  the solution of (\ref{eqna1})-(\ref{eqna2}).  For  small initial mesh-size  $h_1>0$
  there holds
\begin{eqnarray}\label{estimator1}
\norm{{\bf A^{-{1}/{2}}(p-p_h)}}^2+\norm{u-u_h}^2\leq C_{rel} (\eta_h^2+\mu_h^2), 
\end{eqnarray}
where  $0<h\leq h_1$ and $\eta_h$, $\mu_h$  are as defined in (\ref{eta}) and (\ref{mu}).
 \end{theorem}\\
 {\em Proof.} 
 Consider the Helmholtz decomposition: ${\bf e_p}={\mathbf A} \nabla z+\text{Curl} ~\beta$ for $z \in H^1_0(\Omega)$
  and $\beta \in H^1(\Omega)/\mathbb R$. Then
 \begin{equation}\label{ep11}
  (\mathbf{A}^{-1}{\bf e_p},{\bf e_p})=({\bf e_p},\nabla z)+(\mathbf{A}^{-1}{\bf e_p}, \text{Curl}~\beta).
 \end{equation}
 For the  first term  on the right-hand side of (\ref{ep11}), an integration by parts with  (\ref{eq3}) and  the fact ${\rm div~} {\bf p}_h+\gamma u_h=f_h$ lead to
 \begin{align}\label{ep12}
  ({\bf e_p},\nabla z)&=(\text{div}~{\bf e_p},z)=(f-f_h,z)-(\gamma (u-u_h),z)\nonumber\\
 & =(f-f_h,z-\Pi_h z)-(\gamma e_u,z),\nonumber \\
 &\lesssim {osc}_h (f)\;\norm{z}_1+\norm{ e_u} \norm{z}.
 \end{align}
Define $\beta_h:=I_h \beta$, where $I_h:H^1(\Omega)\rightarrow P^1(\mathcal T_h)\cap H^1_0(\Omega)$ is the Clement's interpolation operator \cite{verf}. With
$ {\rm Curl}~ \beta_h \in RT_0(\mathcal T_h)$, ${\rm Curl}~ \beta_h \bot ~\nabla H^1_0(\Omega )~(\bot~ {\rm denotes}~ L^2(\Omega)${ orthogonality}) and (\ref{eq3}),  the  second term 
 on the right-hand side of (\ref{ep11}) can be written as
 \begin{align*}
 (\mathbf{A}^{-1}{\bf e_p},~ \text{Curl}~\beta)
& =-(\mathbf{A}^{-1}{\bf p_h}+ u_h{\bf b}^*,~\text{Curl}~(\beta-\beta_h))-(e_u {\bf b^*},~\text{Curl}~\beta).
 \end{align*}
From the integration by part formula
\begin{equation}\label{ep13}
 (\mathbf{A}^{-1}{\bf e_p},~ \text{Curl}~\beta)
= \sum_{E\in \mc E_h} \int_{ E} [\mathbf{A}^{-1}{\bf p_h}+ u_h{\bf b}^*]\cdot \tau_E (\beta -\beta_h) ds - (e_u {\bf b^*},~\text{Curl}~\beta).
\end{equation}
With the interpolation estimates $\norm{\beta-\beta_h}_E\leq C h_E^{1/2}\norm{\beta}_{1,\omega_E}$, 
    the bounds 
  $\norm{\nabla \beta}_{\omega_E}=\norm{{\rm Curl}~ \beta}_{\omega_E}\leq \norm{\mathbf{A}^{-{1}/{2}} e_p}_{\omega_E}$,  and $\norm{z}\lesssim \norm{z}_1\lesssim \norm{{\bf A^{-1/2} e_p}}$, (\ref{ep11})-(\ref{ep13}) result in
 \begin{align}\label{p_esti}
  \norm{\mathbf{A}^{-{1}/{2}}{\bf e_p}}&\lesssim osc_h(f)
 +\norm{h_E^{1/2}[{\mathbf A}^{-1}{\bf p_h}+u_h{\bf b}^*]\cdot \tau_E}_{\mathcal E_h}+ \norm{e_u}.
 \end{align}
 To estimate   $\norm{e_u},$ start with the triangle inequality 
 \begin{equation}\label{esti_5}
 \norm{e_u}\leq\norm{u-u_h^N}+\norm{u_h^N-u_h}.
 \end{equation}
 For $\tilde e= u_{CR}-u_h^N,$  Lemma 4.5 of \cite{CC} with sufficiently small mesh-size $h$ shows  
 \begin{eqnarray}\label{esti_9}
  \tnorm{\tilde e}_{NC}+\norm{\tilde e} &\lesssim &  osc_h(f).
 \end{eqnarray}
For any $\epsilon >0$, from lemma 3.3 of \cite{CC}, there exists  small mesh-size $h$  such that
\begin{equation}\label{dual}
\norm{u-u_{CR}} \leq \epsilon \tnorm{u-u_{CR}}_{NC} ~~{\rm holds}.
\end{equation}
 A repeated use of the triangle inequality yields estimates for $\norm{u-u_h^N}$ in (\ref{esti_5}) as
\begin{eqnarray}\label{esti_1}
   \norm{u-u_h^N}&\leq& \norm{u-u_{CR}}+\norm{u_{CR}-u_h^N}\nonumber\\
&\lesssim &\epsilon (\tnorm{u-u_h^N}_{NC}+ \tnorm{u_h^N-u_{CR}}_{NC})   +\norm{u_{CR}-u_h^N}\nonumber\\
   &\lesssim & \epsilon \norm{\nabla_{NC}(u-u_h^N)}+ {osc}_h(f).
  \end{eqnarray}
Let   ${\bf p}=-({\bf A}\nabla u+u {\bf b} )$ and define ${\bf \tilde p}=-(\mathbf A\nabla_{NC}u_h^N+u_h{\bf b} )$.  Along with an 
 addition and subtraction of the term  $\mathbf A^{-1} p_h,$ 
\begin{align}\label{esti_6}
 \norm{{\nabla_{NC}( u-u_h^N)}}&\lesssim \norm{\mathbf A^{-1/2}{\bf e_p }} +\norm{e_u}+ \norm{{\bf p_h}-{\bf \tilde p}}.
  \end{align}
 For the third term on the right-hand side of (\ref{esti_6}), (\ref{pmt}) leads to 
 \begin{equation}\label{esti_7}
\norm{{\bf p_h}-{\bf \tilde p}}\leq \norm{(f_h-\gamma u_h)({\bf x}-\text{mid}(T))}\lesssim \norm{h(f_h-\gamma\; u_h)}.
\end{equation}
 The combination of (\ref{esti_1})-(\ref{esti_7}) results in 
 \begin{align}\label{esti_119}
 \norm{u-u_h^N}\lesssim &~ {osc}_h(f)+\epsilon \big(\norm{\mathbf{A^{-1/2}}{\bf e_p}}+ \norm{e_u}\big)+ \epsilon\norm{h(f_h-\gamma\;u_h)}.
 \end{align}
To bound $\norm{u_h^N-u_h}$ in (\ref{esti_5}), use (\ref{umt}), the triangle inequality and the fact that $S(T)\approx h^2$ to obtain
 \begin{eqnarray}\label{esti_3}
  \norm{u_h^N-u_h}&\leq &\left(1+\frac{S(\mathcal{T})}{4} \gamma \right)^{-1} \norm{u_h^N-\Pi_0 u_h^N+\frac{S(\mathcal{T})}{4}(\gamma u_h^N- f_h )}\nonumber\\
 &\lesssim & \norm{h\nabla_{NC}u_h^N} + \tnorm{h^2(u_h^N- u_h)}_{NC}+\norm{h^2\;(f_h-\gamma u_h)}.
 \end{eqnarray}
 A use of (\ref{esti_9})-(\ref{esti_3}) in  (\ref{esti_5}) along with (\ref{eqna2}) leads to 
 \begin{equation}\label{esti_8}
 \norm{e_u}\lesssim {osc}_h (f)+ \epsilon \Big(\norm{\mathbf{A^{-1/2}}{\bf e_p}}+ \norm{e_u}\Big)+ \epsilon\norm{h~{\rm div}~{\bf p_h}}
 +\norm{h(\mathbf{A}^{-1}{\bf p_h}+u_h{\bf b}^*)}.
 \end{equation}
 For  small mesh-size $h_1>0$ with $0 < h\leq h_1,$  (\ref{esti_8}) and (\ref{p_esti}) prove (\ref{estimator1}). 
 \hfill{$\Box$}
\\
\begin{lemma}\label{lmeff}(Efficiency)~Let $({\bf p},u)$ be the solution of (\ref{eq4}) and $({\bf p_h},u_h)$ be the solution of (\ref{eqna1})-(\ref{eqna2})
over the triangulation $\T_h$. Then, it holds
\begin{align}\label{edeff}
C_{\rm eff}(\norm{h_E^{1/2}[{\bf A^{-1}p_h+}u_h{\bf b^*}]\cdot \tau_E}_{\mc E_h}+ &\norm{h({\bf A^{-1}p_h+b^*}u_h)})\nonumber\\
 &\leq  \norm{u-u_h}+\norm{{\bf A^{-1/2}(p-p_h)}}.
 \end{align}
\end{lemma}
{\it Proof.} The proof is divided into two steps.\\
{\em Step 1}. 
Let $b_E$ denote the continuous edge bubble function satisfying $0\leq b_E\leq 1$ on $\omega_E$ and $b_E\in P_2(T)$ for each $T \subset \omega_E$. Let $\phi:= [{\mathbf A}^{-1}{\bf p_h}+u_h{\bf b}^*]\cdot \tau_E $ be a polynomial function along $E$. There exists an extension operator $P:\mathcal C(E)\rightarrow \mathcal C(\omega_E)$ \cite{verf}, 
where $\mathcal C(E)~({\rm resp}. ~\mathcal C (\omega_E))$ denotes the space of the continuous 
functions defined on $E$ (resp. $\omega_E$) such that the  operator $P$ satisfying $P\phi|_E=\phi$ and
\begin{equation}\label{eff_2}
h_E^{1/2}\norm{\phi}_E \lesssim \norm{b_E^{1/2}P\phi}_{\omega_E}\lesssim h_E^{1/2}\norm{\phi}_E.
\end{equation}
An equivalence of norm argument  implies
\begin{equation*}
\norm{\phi}^2_E\lesssim \norm{b_E^{1/2}\phi}_E^2= \int_E (b_EP\phi)[{\mathbf A}^{-1}{\bf p_h}+u_h{\bf b}^*]\cdot \tau_E~ ds.
\end{equation*}
An integration by parts and a use of  ${\rm Curl}({\mathbf A}^{-1}{\bf p_h}+u_h{\bf b}^*)=0 $ result in
\begin{equation}\label{eff_42}
\norm{\phi}^2_E\lesssim \norm{b_E^{1/2}\phi}_E^2\leq -\int_{\omega_E} {\rm Curl}(b_EP\phi)\cdot ({\mathbf A}^{-1}{\bf p_h}+u_h{\bf b}^*)~ ds.
\end{equation}
Note that $(\nabla u, {\rm Curl}(b_E^{1/2}P\phi))_{\omega_E}=0$. Hence, with the help of the 
Cauchy-Schwarz inequality,  (\ref{eff_42}) can be written as
\begin{equation*}
\norm{\phi}^2_E\lesssim \norm{{\mathbf A}^{-1}{\bf p_h}+u_h{\bf b}^*+\nabla u}_{\omega_E} \norm{{\rm Curl}(b_EP\phi)}_{\omega_E}.
\end{equation*}
The inverse inequality, (\ref{eff_2}) 
and a utilization of the definition ${\bf p}=-(\mathbf A \nabla u+u{\bf b} )$, result in
\begin{align*}
\norm{h_E^{1/2}[{\bf A^{-1}p_h}+u_h{\bf b^*}]\cdot \tau_E}_{E}~~\lesssim \norm{ {\bf p-p_h}}_{\omega_E} +\norm{u-u_h}_{\omega_E}.
\end{align*}
A summation over all the edges leads to an estimate of the first term on the left-hand side of (\ref{edeff}).\\
{\em Step 2}. Define the function ${\textbf q_T}:= b_T(\mathbf A^{-1}{ \bf p_h} +u_h{\textbf b^*} )\in P_4(T)\cap W^{1,\infty}_0(T)$ and the cubic bubble function 
$
 ~b_T = 27\lambda_1 \lambda_2 \lambda_3 \in P_3(T)\cap C_0(T)
$
in terms of the barycentric coordinates $\lambda_1, \lambda_2, \lambda_3$ of $T \in \mathcal T_h$ \cite{verf}.
Since ${\mathbf A}^{-1}{\bf p_h}+u_h{\bf b}^*$ is affine on $T\in \mathcal T_h$, an equivalence of norm argument shows
\begin{align*}
 \norm{\mathbf A^{-1}{ \bf p_h} +u_h{\textbf b^*} }_T^2&\lesssim  \int_T{\textbf q_T}\cdot({ \mathbf A^{-1}\bf p_h} +u_h{\textbf b^*} )~dx.
 \end{align*}
The definition of ${\bf p}$ and (\ref{eq3}) show that 
\begin{align*}
\norm{{\mathbf A^{-1} \bf p_h} +u_h{\textbf b^*} }_T^2 
& \lesssim \int_T{\textbf q_T}\cdot \big ( \mathbf A^{-1}({ \bf p_h - p}) -(u-u_h){\textbf b^*} \big)~dx
~-\int_T{\textbf q_T}\cdot \nabla u ~dx.                      
 \end{align*}
The Cauchy-Schwarz inequality with $\norm{{\textbf q_T}}_T\lesssim \norm{\mathbf A^{-1}_h{\bf p_h}+u_h{\textbf b_h^*} }_T$ is employed in the first 
two terms. Adding the zero terms $\nabla u_h|_T$ to the right-hand side of the above equation,  an integration by parts  
 shows  that 
\begin{eqnarray*}
 h_T\norm{{\mathbf A^{-1} }{\bf p_h}+u_h
 {\bf b} }^2_T &\lesssim &h_T\norm{\mathbf A^{-1}{\bf p_h}+u_h{\bf b}^* }_T\Big ( \norm{{\bf p-p_h}}_T
+\norm{u-u_h}_T \Big )\nonumber\\
 &&+ h_T \int_T (u-u_h)\text{div}~{\bf q}_T ~dx.
 \end{eqnarray*}
Since $\textbf q_T \in P_4(T)$, an inverse estimate yields
\begin{equation*}
 h_T\norm{\text{div}~ \textbf q_T}_T\lesssim \norm{\textbf q_T}_T\lesssim \norm{\mathbf A^{-1}{\bf p_h}+u_h{\bf b}^* }_T.
\end{equation*}
Since $ h_T \lesssim 1$, it follows
\begin{align*}
h_T\norm{\mathbf A^{-1}{\bf p_h}+u_h{\bf b}^* }_T&\lesssim \norm{u-u_h}_T+\norm{{\bf p}-{\bf p_h}}_T.
\end{align*}
A summation  over all elements leads to an estimate of  second  term on the left-hand side of (\ref{edeff}). This
concludes the proof. \hfill{$\Box$}
\\
\\
 Let $\mc T_H$ and $\mc T_h$ with $H<h$ denote  nested triangulations, 
  and $({\bf p_h}, u_h)$ and $({\bf p_H}, u_H)$ denote the solutions of  (\ref{eqna1})-(\ref{eqna2}) obtained with right-hand sides $f_h$ and $f_H$ over
  $\mathcal T_h$ and $\mathcal T_H$, respectively.
 The following notations are used  in the sequel: 
\begin{subequations}\label{com_1}
\begin{eqnarray}
&&E_H^2:=\norm{{\bf A}^{-1/2}{\bf (p_h-p_H)}}^2+\norm{u_h-u_H}^2,\\
&&e_{p_h}^2:=\norm{{\bf A}^{-1/2}{\bf  (p-p_h)}}^2,\qquad\qquad e_{p_H}^2:=\norm{{\bf A}^{-1/2}{\bf  (p-p_H)}}^2,\\
&&e_{u_h}^2:=\norm{u-u_h}^2, \qquad~~~~~~~~~~~\qquad e_{u_H}^2:=\norm{u-u_H}^2,\\
&&e_h^2~:=e_{p_h}^2+e_{u_h}^2,\qquad\qquad ~~~~~~~~~~~~~ e_H^2~:=e_{p_H}^2+e_{u_H}^2.
\end{eqnarray}
\end{subequations}

\begin{lemma}\label{lm3}(Volume estimator reduction) For given $0<\theta_B \leq 1$, there exist constants $0<\delta_1,~ \rho_B < 1$ and  the positive constant $\varLambda_2$ such that $\mu_H$ and $E_H$ defined in (\ref{mu}) and (\ref{com_1}) satisfy
\begin{eqnarray}
&&  \mu_h^2 \leq (1+\delta_1) \mu_H^2 +\varLambda_1 H^2 E_H^2 \qquad \qquad {\rm for~~the~~ {\bf Case ~(A)}},\label{et_rd_11} \\
&&  \mu_h^2 \leq \rho_B \mu_H^2 +\varLambda_1 H^2 E_H^2 \qquad \qquad \qquad {\rm for~~the ~~{\bf Case ~(B)}}.\label{et_rd11} 
\end{eqnarray}
\end{lemma}
{\it Proof}.  
For any  triangle $K \in \mathcal T_H$, which gets refined in the level $\mc T_h$, that is, $K\in \mc T_H\setminus \mc T_h$, there exist triangles  $T_1,T_2,\cdots,T_J $
 such that $K= T_1\cup T_2 \cdots \cup T_J$. For $K \in \mathcal T_H \setminus \mathcal T_h$ 
 \begin{eqnarray*}
\mu^2_h(K)&=&\sum_{j=1}^J\Big( \norm{h(f-f_h)}_{T_j}^2+\norm{h~{\rm div}~{\bf p_h}}^2_{T_j}+\norm{h ({\bf A^{-1}p_h }+u_h{\bf b^*})}^2_{T_j}\Big).
\end{eqnarray*}
At least one refinement of $T$ implies $h\leq H/2$,  the fact $\norm{f-f_h}_{T_j}\leq\norm{f-f_H}_{T_j}$  along with the triangle 
inequality yields
\begin{eqnarray}
\mu^2_h(K)&\leq& \frac{H^2}{4}\sum_{j=1}^J\Big( \norm{f-f_H}_{T_j}^2+\norm{~{\rm div}~{(\bf p_h- p_H)}+{\rm div}~{\bf p_H}}^2_{T_j}\nonumber \\
&&+\norm{({\bf A^{-1}p_H }+u_H{\bf b^*})+{\bf A^{-1}(p_h-p_H) }+(u_h-u_H){\bf b^*}}^2_{T_j}\Big)\nonumber
\end{eqnarray}
For each $K\in \T_H$, a use of $\norm{f_h-f_H}_K\leq \norm{\Pi_h(1-\Pi_H)f}_K\leq \norm {f-f_H}_K$ and  (\ref{eqna2}), with  Young's inequality yields for $\delta_1>0$
\begin{align}\label{51}
\mu^2_h(K)&\leq \frac{H^2}{4}\sum_{j=1}^J\big( \norm{(f-f_H)}_{T_j}^2+2\norm{f_h-f_H}^2_{T_j}+ 3\norm{{\rm div}~{\bf p_H}}^2_{T_j}+6\norm{\gamma (u_h-u_H)}^2\nonumber \\
&~+(1+\delta_1)\norm{{\bf A^{-1}p_H }+u_H{\bf b^*}}^2_{T_j}+(1+\frac{1}{\delta_1})\norm{{\bf A^{-1}(p_h-p_H) }+(u_h-u_H){\bf b^*}}^2_{T_j}\big)\nonumber\\
&\leq \frac{3}{4}(1+\delta_1) \mu_H^2(K) +(1+\frac{1}{\delta_1}) C_1 H^2 E_H^2(K),
\end{align} 
where $C_1=\max\{\norm{{\bf A}^{-1/2}}_\infty^2,\norm{{\bf b^*}}_\infty^2+6 \norm{\gamma}_\infty^2\}$. Denote  $(1+\frac{1}{\delta_1}) C_1=:\varLambda_1.$\\
For $K\in \mathcal T_H \cap \mathcal T_h$,
 \begin{eqnarray}\label{52}
\mu^2_h(K)
&\leq& (1+\delta_1)\mu_H^2(K) +\varLambda_1 H^2 E_{H}^2{(K)}.
 \end{eqnarray}
From (\ref{51}) and (\ref{52}), a summation over all triangles implies 
 \begin{eqnarray}\label{et_rd1}
\mu_h^2 \leq (1+\delta_1) \mu_H^2+\varLambda_1  H^2 E_H^2,
\end{eqnarray}
for  the  {\bf Case (A)}.\\
For  the {\bf Case (B)}, a summation over all triangles, 
the marking criteria for $\mc M_H\subset \T_H$, 
\begin{equation}\label{81}
\theta_B \mu_H^2 \leq \mu_H^2(\mathcal M_H) \leq \mu_H^2 (\mathcal T_H \setminus T_h).
\end{equation}
  A use of (\ref{51}), (\ref{52}) and (\ref{81})  leads to the sharper bound 
\begin{eqnarray*}\label{vol4}
\mu_h^2 & =&\mu_h^2(\mathcal T_H \setminus T_h) +\mu_h^2(\mathcal T_h \cap \mathcal T_H)\nonumber \\
&\leq& \frac{3}{4}(1+\delta_1) \mu_H^2(\mathcal T_H \setminus T_h) +(1+\delta_1)\mu_H^2(\mathcal T_h \cap \mathcal T_H)+\varLambda_1 H^2 E_H^2.\nonumber\\
&\leq &(1+\delta_1)\mu_H^2- \frac{1}{4}(1+\delta_1) \mu_H^2(\mathcal T_H \setminus T_h)+\varLambda_1  H^2 E_H^2\nonumber\\
 & \leq &(1+\delta_1)(1-\frac{\theta_B}{4})\mu_H^2+\varLambda_1  H^2 E_H^2.
\end{eqnarray*}
For given  $0<\theta_B<1$, the selection  $0<\delta_1< {\theta_B}/{(4-\theta_B)}$  results in  \\
$0<\rho_B=(1+\delta_1)(1-\frac{\theta_B}{4})<1$.~
 This concludes the proof.\hfill${\Box}$ 
\begin{remark}\label{rmk3} Irrespective of the marking criteria, from (\ref{et_rd_11})-(\ref{et_rd11}), it
 follows that 
$$\mu_h^2 \leq 2 \mu_H^2+\varLambda_1  H^2 E_H^2.$$
\end{remark}
\begin{corollary}\label{lm13} Let  $\T_h$ be a refined triangulation of $\T_H$. Then it holds,
 \begin{eqnarray}\label{lm12}
\frac{\gamma_0^2}{2}\mu^2_H\leq \mu_h^2  +2 \norm{H(f_h-f_H)}^2+C_1 H^2 E_H^2,
 \end{eqnarray}
 where for $T\in \T_h$,  $T' \in \T_H$ and $T\subset \T'$, $\gamma_0$ is defined as
 $0<\gamma_0<1$ if $T'$ get refined, otherwise $\gamma_0=1$.
\end{corollary}

\begin{lemma}\label{lm2}( Edge estimator reduction) Given $0<\theta_A \leq 1$, there exist constants $0<\delta_2,~ \rho_A < 1$  and  a positive constant $\varLambda_2$ such  $\eta_H$ and $E_H$ defined in (\ref{eta}) and (\ref{com_1})
satisfy
\begin{eqnarray}
  \eta_h^2 \leq \rho_A \eta_H^2 +\varLambda_2 E_H^2 \qquad \qquad {\rm for~~the~~{\bf Case~ (A)}},\label{et_rd2} \\  
    \eta_h^2 \leq (1+\delta_2)\eta_H^2 +\varLambda_2 E_H^2 \qquad \qquad ~ {\rm for~~the ~~{\bf Case ~(B)}}.\label{et_rd_2}
\end{eqnarray}
\end{lemma}
{\it Proof}. For all $F \in \mathcal E_H$, either $F\in \mathcal E_h$ or there exist $E_1, E_2, ... E_J \in \mathcal E_{h}$  with 
$F = E_1\cup E_2\cup ...\cup E_J$ for $J \geq 2$. For the case  $F \in \mathcal E_H \setminus \mathcal E_h$, a use of Young's inequality yields 
\begin{eqnarray}\label{99}
\eta_h^2(F)&=&\sum_{j=1}^J \eta^2_h (E_{j})=\sum_{j=1}^J \norm{h_{E_j}^{1/2}[{\bf A^{-1}p_h}+u_h{\bf b^*}]\cdot \tau_{E_j}}^2_{E_j} \nonumber\\
&\leq & (1+\delta_2)\frac{H_{F}}{2}\norm{[{\bf A^{-1}p_H}+u_H{\bf b^*}]\cdot \tau_{E}}^2_F \nonumber\\
&&+(1+\frac{1}{\delta_2})\sum_{j=1}^J {h_{E_j}}\norm{[{\bf A^{-1}(p_h-p_H)}+(u_h-u_H){\bf b^*}]\cdot \tau_{E_j}}^2_{E_j}.
\end{eqnarray}
Here $H_F$ denotes the lenght of edge $F$. The fact that  $J\geq 2$ implies that there exists  at least one bisection of $F \in \mathcal E_H \setminus \mathcal E_h$.\\
For $F\in \mc E_H \cap \mc E_h$,
\begin{eqnarray}\label{98}
\eta_h^2(F)
&\leq & (1+\delta_2)\norm{{H_{F}}[{\bf A^{-1}p_H}+u_H{\bf b^*}]\cdot \tau_{E}}^2_E \nonumber\\
&&+(1+\frac{1}{\delta_2}) {H_{F}}\norm{[{\bf A^{-1}(p_h-p_H)}+(u_h-u_H){\bf b^*}]\cdot \tau_{F}}^2_{F}.
\end{eqnarray}
For any $E\in \mathcal E_h$ with $E \nsubseteq \cup \mathcal E_H$,  $E$ is the interior edge of some element of $T\in \T_H$ and hence,  $ [{\bf A^{-1}p_H}+u_H{\bf b^*}]\cdot \tau_E =0$ along $E$ implies  $\eta_h^2(E)=\norm{h_E[{\bf A^{-1}(p_h-p_H)}+(u_h-u_H){\bf b^*}]\cdot \tau_E}^2_E$. 
Now, consider this and the cases (\ref{99})-(\ref{98})   to obtain
\begin{eqnarray*}
\eta_h^2 =\sum_{E \in \mathcal E_h} \eta_h^2(E) &\leq &\frac{1+\delta_2}{2} \sum_{E\in \mathcal E_H \setminus \mathcal E_h} \eta_H^2(E)
+ (1+\delta_2) \sum_{E\in \mathcal E_H \cap \mathcal E_h} \eta_H^2(E)\nonumber \\ &&~+ (1+\frac{1}{\delta_2})\sum_{E \in \mathcal E_h} 
\norm{h_E [{\bf A^{-1}(p_h-p_H)}+(u_h-u_H){\bf b^*}]\cdot \tau_E}^2_{E}.
\end{eqnarray*}
The inverse inequality $\norm{{\bf q_H}}^2_{{E}} \lesssim H^{-\frac{1}{2}}\norm{{\bf q_H}}^2_{{\omega_E}}$ results in
\begin{eqnarray*}
\norm{h_E [{\bf A^{-1}(p_h-p_H)}+(u_h-u_H){\bf b^*}]\cdot \tau_E}^2_{E}
&\lesssim &  \norm{ \big({\bf A^{-1}(p_h-p_H)}+(u_h-u_H){\bf b^*}\big)}^2_{\omega_E}
\end{eqnarray*}
for the edge patch $\omega_E$ of $E$ in $\mathcal T_h$. Since there is only a finite overlap of all edge patches, there holds
\begin{eqnarray*}
 \eta_h^2\leq \frac{1+\delta_2}{2}\sum_{E \in \mathcal E_H \setminus \mathcal E_h } \eta_H^2+(1+\delta_2)\sum_{E \in \mathcal E_H \cap \mathcal E_h } \eta^2_H
+C_2 (1+\frac{1}{\delta_2}) E_H^2.
\end{eqnarray*}
Denote $C_2 (1+\frac{1}{\delta_2})=:\varLambda_2$ to obtain 
\begin{eqnarray}\label{edge_4}
 \eta_h^2\leq \frac{1+\delta_2}{2}\sum_{E \in \mathcal E_H \setminus \mathcal E_h } \eta_H^2+(1+\delta_2)\sum_{E \in \mathcal E_H \cap \mathcal E_h } \eta^2_H
+\varLambda_2 E_H^2.
\end{eqnarray}
 For the {\bf Case (A)}, the marking criteria  for $\mc M_H\in \mc E_H$,   $$\theta_A \eta_H^2 \leq \eta_H (\mathcal M_H)\leq \eta (\mathcal E_H \setminus \mathcal E_h)$$
 leads to 
\begin{eqnarray*}
\eta_h^2 &\leq & (1+\delta_2) \eta_H^2 -\frac{(1+\delta_2)}{2} \sum_{E \in \mathcal E_H \setminus \mathcal E_h}\eta_H^2+\varLambda_2E_H^2 \nonumber\\
& \leq &(1+\delta_2)(1-\frac{\theta_A}{2}) \eta_H^2 + \varLambda_2 E_H^2.
\end{eqnarray*}
For any given  $\theta_A$, the choice of $\delta_2={\theta_A}/{(4-2\theta_A)}$ implies $0<\rho_A=(1+\delta_2)(1-\frac{\theta_A}{4})<1$.  Hence, (\ref{et_rd2}) holds. 
For  the {\bf Case (B)},  a summation over the all edges  implies 
\begin{equation*}\label{et_rd3}
 \eta_h^2 \leq (1+\delta_2) \eta_H^2 + \varLambda_2 E_H^2.
\end{equation*}
This completes the rest of the proof. \hfill${\Box}$ 
\begin{lemma}\label{lm6}(Quasi-orthogonality)
Let $\mathcal T_h$  be  a refined triangulation of $\mathcal T_H$. Then for  small  initial mesh-size $h\leq h_2$, there exist constants $0<\alpha_1 ,\alpha_3<1$ such that
\begin{eqnarray}\label{er_re}
(1-\alpha_1){ e_{h}}^2\leq { e_H}^2-\alpha_3 E_H^2+ \varLambda_3 \mu_H^2.
\end{eqnarray}
\end{lemma}
{\it Proof}. The following hold 
\begin{align}
~\norm{{\bf A^{-1/2} (p- p_h)}}^2&=\norm{{\bf A^{-1/2} (p- p_H)}}^2-\norm{{\bf A^{-1/2} (p_h- p_H)}}^2\nonumber \\
&~~-2\big({{\bf A^{-1} (p- p_h), p_h-p_H}}\big),\label{nov2}\\
\qquad\qquad \norm{{u- u_h}}^2&=\norm{{u- u_H}}^2-\norm{{ u_h- u_H}}^2-2({{ u- u_h, u_h-u_H}}).\label{nov3}
 \end{align}
The term $\big({{\bf A^{-1} (p- p_h), p_h-p_H}})$ is estimated  by introducing the intermediate term ~${\bf \tilde p}= -({\bf A}\nabla_{NC}u_h^N+u_h{\bf b})$ as  
 \begin{eqnarray}\label{ortho1}
({{\bf  A^{-1} (p- p_h), p_h-p_H}}) &=& ({{\bf A^{-1} (p- \tilde p), p_h-p_H}})+({{\bf A^{-1} (\tilde p- p_h), p_h-p_H}})\nonumber\\
&=& -(\nabla_{NC} (u-u_h^N), {\bf p_h- p_H})-((u-u_h){\bf b^*},{\bf p_h-p_H})\nonumber\\
&&+({{\bf A^{-1} (\tilde p- p_h), p_h-p_H}}).
 \end{eqnarray}
An elementwise integration by parts of the first term on the  right-hand side of (\ref{ortho1})  yields
\begin{align}\label{ortho2}
-(\nabla_{NC}(u-u_h^N), {\bf p_h- p_H})&=(u-u_h^N, \text{div}~ ({\bf p_h- p_H}))\nonumber\\
&~~~-\sum_{E\in \mathcal E_h} \int_E [u-u_h^N]({\bf p_h -p_H})\cdot \nu_E ~ds.
\end{align}
The second term on right-hand side of (\ref{ortho2}) is zero, as $({\bf p_h -p_H})\cdot \nu_E$ is continuous along the edge $E$ and constant on $E$, and  $u_n^N\in CR_0^1 (\T)$, $\int_E [u-u_h^N] ds =0$. 
The second equation of weak formulation (\ref{eqna2}) and the definition of $L^2$-projection show
\begin{align}\label{ortho3}
 (u-u_h^N, & \text{div}~ ({\bf p_h- p_H})) \nonumber \\
 &= ((u-u_h^N)-\Pi_H (u-u_h^N), f_h-f_H)-(u-u_h^N,\gamma (u_h-u_H)).
\end{align}
A substitution of (\ref{ortho2})-(\ref{ortho3}) in (\ref{ortho1}) with a use of   the Cauchy-Schwarz inequality and the Poincar$\acute{e}$ inequality results in 
\begin{align*}
&({{\bf  A^{-1} (p- p_h), p_h-p_H}}) \lesssim   \norm{\nabla_{NC} (u-u_h^N)} \norm{H(f_h-f_H)}+\norm{u-u_h^N} \norm{u_h-u_H}  \nonumber\\
&\hspace{10mm}+\norm{u-u_h} \norm{{\bf A^{-1/2} (p_h-p_H)}}+\norm{{\bf \tilde p- p_h}} \norm{{\bf A^{-1/2} (p_h-p_H)}}.
\end{align*}
Using the estimates (\ref{esti_1})-(\ref{esti_8}), Lemma \ref{lm8} and  the addition of the term
\\ $(u-u_h,u_h-u_H)$ yield
\begin{eqnarray*}
&& 2|({{\bf A^{-1} (p- p_h), p_h-p_H}})|+2|(u-u_h,u_h-u_H)|\leq C_3 (e_{p_h}+e_{u_h})osc_H (f)\nonumber\\
&&~~+C_3(\norm{h f_h}+\norm{h( {\bf A^{-1} p_h}+ u_h{\bf b^*})})osc_H (f) +C_3\epsilon (e_{p_h}+e_{u_h})(\norm{u_h-u_H}\nonumber\\ 
&&~~+\norm{{\bf A^{-1/2}(  p_h-p_H)}})+C_3(osc_h (f)+\norm{hf_h}+\norm{h{\bf A^{-1} p_h}+ u_h{\bf b^*}})\nonumber \\
&&~~(\norm{u_h-u_H}+\norm{\bf A^{-1/2}(  p_h-p_H)}).
\end{eqnarray*}
The Young's inequality and rearrangement of terms result in 
\begin{eqnarray}\label{pp1}
 &&2|({{\bf A^{-1} (p- p_h), p_h-p_H}})|+2|(u-u_h,u_h-u_H)|  \nonumber\\
&&\qquad\qquad\leq (\delta_3+{4C_3^2\epsilon^2})(e_{p_h}^2+e_{u_h}^2)+\frac{1}{2} E_H^2
+7C_3^2 \mu_h^2+\big(\frac{1}{2}+\frac{C_3^2}{2\delta_3}\big)osc_H ^2(f).
 \end{eqnarray}
As $\mathcal T_h$ and $\mathcal T_H$ are the nested triangulations, Remark \ref{rmk3}  implies that \\ $\mu_h^2 \leq  2\mu_H^2 +\varLambda_1 H^2 E_H^2$. 
Hence, with notations  $0<\alpha_1=(\delta_3+{4C_3^2\epsilon^2}),$ $\alpha_2= (\frac{1}{2}+{7C_3^2}\varLambda_1 H^2)$  and $\varLambda_3 =({14C_3^2}+\frac{1}{2}+\frac{C_3^2}{2\delta_3})$, (\ref{pp1})  reduces to
\begin{align}\label{nov1}
 2|({{\bf A^{-1} (p- p_h), p_h-p_H}})|+& 2|(u-u_h,u_h-u_H)| \nonumber \\
 & \leq \alpha_1(e_{p_h}^2+e_{u_h}^2)+ \alpha_2 E_H^2+\varLambda_3\mu^2_H.
\end{align}
A combination of (\ref{nov2})-(\ref{nov3}) with (\ref{nov1}) yields
\begin{eqnarray}
({ e_{p_h}}^2+e_{u_h}^2)\leq ({ e_{p_H}}^2+e_{u_H}^2)- E_H^2+\alpha_1(e_{p_h}^2+e_{u_h}^2)+ \alpha_2 E_H^2+\varLambda_3\mu^2_H.\nonumber
 \end{eqnarray}
Hence $ (1-\alpha_1)e_{h}^2\leq e_{H}^2-(1-\alpha_2) E_H^2+ \varLambda_3 \mu_H^2.$\\
 For any $\epsilon >0$,  
it always possible to find  a  small initial mesh-size  $h_2$ and $0<\delta_3<1$, such that
$0<\alpha_1=\delta_3+4C_3^2\epsilon^2<1$ and $0<\alpha_3=1-\alpha_2=1/2-7C_3^2\varLambda_1 H^2$. 
Hence, (\ref{er_re}) holds true and this completes the rest of the proof.\hfill${\Box}$ 
%
\begin{remark} \label{rmk2} For proving the contraction property, a more sharper bound for \\
$\alpha_1 <\alpha_4:={1}/{2}\min\{1, {\alpha_3(1-\rho_A)}/({(\varLambda_1+\varLambda_2)C_{rel}})\}<1$ will be selected. 
\end{remark}
 \begin{corollary}\label{coro3} Under the assumption that  small initial mesh-size $h_2 >0$, and   the constants  defined in Lemma \ref{lm6}, the following result holds for $0<h\leq h_2$
\begin{eqnarray}
 {\bf e_{H}}^2&\leq & (1+\alpha_1){\bf e_h}^2+(1+\alpha_2) E_H^2+ \varLambda_3 \mu_H^2\nonumber \\
 &\leq&  2{\bf e_h}^2+ 2 E_H^2+ \varLambda_3 \mu_H^2.
\end{eqnarray}
 \end{corollary}
\begin{lemma}({Quasi-discrete reliability})\label{coro2}
~Let $({\bf p}_h, u_h )$ and $({\bf p}_H, u_H)$ be the MFEM solutions of (\ref{eqna1})-(\ref{eqna2}) over the triangulations $\T_h$ and $\T_H$, respectively.  There exists a constant $C_4 >0$ such that for any $\epsilon >0$ 
 \begin{align*}
\norm{{\bf A^{-1/2} (p_h -p_H)}}^2+\norm{u_h-u_H}^2 &\leq C_4\big( \eta_H^2(\mathcal E_H \setminus \mathcal E_h) + \norm{H (f_h- f_H)}^2\\
&~~~ +  \epsilon^2 ( e_h^2+  e_H^2) +\mu_h^2 +\mu_H^2\big),
\end{align*}
where $E_H, e_h, e_H$ are defined  in (\ref{com_1}) and $\mu_h$ in (\ref{mu}).
%
\end{lemma}\\
{\it Proof}.
Introduce the discrete mixed finite element problem: seek
 $ (\tilde{\bf p}_h, \tilde u_h) \in RT_0(\mathcal {T}_h)\times P_0(\mathcal T_h)$ such that
  \begin{eqnarray}
  &&(\mathbf{A}^{-1}\tilde {\bf p}_h,{\bf q_h})-(\text{div}~{\bf q_h},\tilde u_h)=-( u_H {\bf b}^*,{\bf q_h}) \qquad 
  \text{for all}\, {\bf q_h}\in RT_0(\mathcal {T}_h), \label{eqna8}\\
&&  (\text{div}~ \tilde{\bf p}_h,v_h)=(f_H,v_h)- (\gamma u_H, v_h) \qquad  \text {for all} \, v_h \in P_0(\mathcal{T}_h). \label{eqna9}
 \end{eqnarray}
 Define  the nonconforming discrete  problem  corresponding to (\ref{eqna8})-(\ref{eqna9}): seek  $\tilde u_h^N\in CR_0^1(\mathcal T_h)$ as the solution of
 \begin{equation}
  (\mathbf{A}\nabla_{NC} \tilde u_h^N+{u}_H {\bf b},\nabla_{NC} v_{CR})+(\gamma {u}_H,v_{CR})=
  (f_H,v_{CR}), {\forall} v_{CR}\in CR^1_0(\mathcal T_h).\label{eqna10}
 \end{equation}
The solution ${\tilde {\bf p}_h}$ of (\ref{eqna8})-(\ref{eqna9}) can be written in the terms of $\tilde u_h^N$ as
 \begin{equation} \label{pmt1}
  \tilde {\bf p}_h({\bf x})=-\left(\mathbf{A}\nabla _{NC} \tilde u_h^N+u_H {\bf b} \right)+\left(f_H-\gamma  u_H\right)
  \frac{\left({\bf x}-\text{\rm mid}(T)\right)}{2} \text{ for}~~ {\bf x}\in T.
 \end{equation}
 Now for the estimates of $\norm{{\bf p_h- p_H }}$,  use ${\tilde {\bf p}_h}$ as an intermediate term to split  ${\bf p_h}-{\bf p_H}:= ({\bf p_h}-{\bf {\tilde p_h}})+({\bf {\tilde p_h}}-{\bf p_H})$ and  the triangle 
 inequality. 
  From the representation formula (\ref{pmt}) and (\ref{pmt1}) of ${\bf p_h}$ and ${\bf \tilde p_h }$, respectively, it follows that
\begin{eqnarray*}
{{\bf p_h- \tilde p_h }}= - ({\bf A} \nabla (u_{h}^N-\tilde u_h^N)+ (u_h -u_H){\bf b})+\frac{1}{2} (f_h-f_H-
\gamma (u_h-u_H))(x- {\rm mid} (T)).
 \end{eqnarray*}
 The triangle inequality shows 
\begin{align}\label{rea_4}
\norm{  {\bf p_h- \tilde p_h } } \lesssim \norm{{\bf A} \nabla (u_{h}^N-\tilde u_{h}^N)}+\norm{ u_h -u_H}+\norm{h(f_h-f_H)}. 
\end{align}
Subtracting (\ref{eqna3}) from (\ref{eqna10}) leads to
 \begin{align}\label{rea_6}
({\bf A} \nabla (u_{h}^N-\tilde u_{h}^N), \nabla v_{CR})&= (f_h-f_H,v_{CR})+(\gamma(u_h-u_H),v_{CR})\nonumber \\
&~~-((u_h-u_H){\bf b},\nabla v_{CR}).
 \end{align}
A substitution $v_{CR}={\bf A} (u_{h}^N-\tilde u_{h}^N)$ in (\ref{rea_6}) with the $L^2$-projection property and the Cauchy-Schwarz inequality results in
$$\norm{{\bf A}\nabla (u_{h}^N-\tilde u_{h}^N)}\lesssim \norm{H(f_h-f_H)}+\norm{u_h-u_H}.$$ 
With this estimate, (\ref{rea_4}) reduces to 
 \begin{eqnarray}\label{rea_9}
\norm{  {\bf p_h- \tilde p_h } } \lesssim\norm{H(f_h-f_H)}+ \norm{ u_h -u_H}.\label{rea_7}
\end{eqnarray}
For a bound of  the term $ (u_h-u_H),$ a use of (\ref{esti_8}) shows
 \begin{eqnarray}\label{rea_8}
\norm{u_h-u_H}^2 &\lesssim& \norm{u_h-u}^2+ \norm {u-u_H}^2 \nonumber \\
&\lesssim & \epsilon^2 e_h^2 +\mu_h^2 + \epsilon^2 e_H^2+\mu_H^2.
 \end{eqnarray}
For the estimate of $\norm{{\bf p_h- p_H }}$, note that
${\rm div}~ {\bf \tilde p_h} =f_H-\gamma u_H={\rm div}~ {\bf p_H}$. It implies $\rm {div} ~(\tilde {\bf p}_h-{\bf p}_H)=0$, and hence,  
$ {\bf \tilde p_h- p_H}$ is a piecewise constant vector function over $\mathcal T_h$.
  The discrete Helmholtz decomposition  states for 
 ${\bf A (\tilde p_h-p_H)= A}\nabla \alpha_{CR}+ \rm {Curl}~\beta_h$, where $\alpha_{CR}\in CR_0^1(\mathcal T_h)$ and $\beta_h \in P_1(\mathcal T_h)\cap C(\bar \Omega)$, and hence, 
   \begin{eqnarray*}
 \norm{{\bf \tilde p_h -p_H}}^2=({\bf \tilde p_h -p_H}, \nabla \alpha_{CR}+ A^{-1} \rm {Curl}~ \beta_h ).
 \end{eqnarray*}
 Let $\beta_H := I_H \beta_h$ be the Scott-Zhang quasi-interpolation operator  for  $E \in \mc E_H$, 
 where $\mathcal E_H$ the set of  edges on the triangulation $\mathcal T_H$ and its neighbourhood $\omega_E$ with
  \begin{equation}\label{drel22}
  \norm{\beta_h- \beta_H}_E\leq C h_E^{{1}/{2}} \norm{\beta_h}_{H^1{(\omega_E})}.
  \end{equation}
 Note that $\norm{\beta_h- \beta_H}_E=0 $ if $E \in \mc E_h\cap\mc E_H.$\\
 The fact $(({\bf \tilde p_h -p_H}), \nabla_{NC} \alpha_{CR})=0$ shows~$ \norm{{\bf \tilde p_h -p_H}}^2=({\bf \tilde p_h -p_H, A^{-1}}\rm {Curl}~ \beta_h )$. 
The weak formulation (\ref{eqna1}) with  ${\bf q_H}= {\rm Curl}~ \beta_H \in RT_0(\mc T_H)\subset RT_0(\mc T_h)$ over $\mc T_H$ and (\ref{eqna8})  with ${\bf q_h}= {\rm Curl}~ \beta_H \in RT_0(\mc T_h)$ and integration by parts lead to
\begin{eqnarray}\label{rea_2}
 \norm{{\bf \tilde p_h -p_H}}^2& = & -( u_H{\bf b^*}, \rm {Curl}~ \beta_h) - ({\bf A^{-1}p_H}, \rm {Curl}~\beta_h)\nonumber \\
 &= &({\bf A^{-1} p_H}+ u_H {\bf b^*}, {\rm Curl}~ (\beta_H-\beta_h)) \nonumber\\
& = & \sum_{E \in \mc E_h} \int_E [{\bf A^{-1}p_H+b^*}u_H]\cdot \tau_E (\beta_H-\beta_h) ~ds \nonumber\\
&&~-\sum_{T \in \mc T_h}\int_T  {\rm Curl}~ ({\bf A^{-1}p_H}+u_H{\bf b^*}) (\beta_H-\beta_h)~ ds. \nonumber
\end{eqnarray}
Since $\rm {Curl}~ ({\bf A^{-1}p_H}+u_H{\bf b^*})=0$
\begin{eqnarray}
\norm{{\bf \tilde p_h -p_H}}^2& \lesssim & \sum_{E \in \mc E_h \setminus \mc E_H } \norm{[{\bf A^{-1}p_H}+u_H {\bf b^*}]\cdot \tau_E}_{E} \norm{\beta_h-\beta_H}_{E}.
\end{eqnarray}
The estimates of $\norm{\beta_h-\beta_H}_{E}$ from (\ref{drel22}) and the bound  $\norm{\nabla \beta_h}=\norm{{\rm Curl}~ \beta_h} \lesssim \norm{{\bf \tilde p_h -p_H}}$ together with (\ref{rea_2}) result in
\begin{equation}\label{rea_1}
 \norm{{\bf \tilde p_h -p_H}} \lesssim \eta_H (\mc E_H \setminus \mc E_h).
\end{equation}
 A combination of (\ref{rea_1}), (\ref{rea_9}) and (\ref{rea_8}) leads to
 \begin{eqnarray}\label{53}
 \norm{  {\bf p_h-  p_H } } ^2\lesssim \eta_H^2(\mathcal E_H \setminus \mathcal E_h) + \norm{H (f_h- f_H)}^2 + \epsilon^2 ( e_h^2+  e_H^2) +\mu_h^2 +\mu_H^2.
 \end{eqnarray} 
 The combination of (\ref{rea_8}) and (\ref{53}) completes the proof.\hfill${\Box}$
\section{Convergence Analysis}
This section is devoted to the convergence analysis of the adaptive mixed finite element method.
\subsection{Contraction property}
Denote two consecutive adaptive loop levels as $\ell$ and $\ell +1$. Let $e_\ell$ , $\eta_\ell$ and $\mu_\ell$  denote
the error and the estimator terms on the level  $\ell$ with triangulation $\mathcal T_\ell$. 
 Based on the reduction properties of the error, the error estimators and quasi-orthogonal property developed in last section, the  contraction property is proved  for 
 the weighted term $\xi^2_\ell$, which is a linear combination of error $e^2_\ell$ and the  estimator terms $\eta^2_\ell$ and $\mu^2_\ell$, between two consecutive adaptive loops.  
 \begin{theorem}\label{thm1}(Contraction Property) Let $\mathcal T_{\ell+1}$ be a refinement of $\mathcal T_{\ell}$ using  AMFEM algorithm. Given $0<\theta_A, \theta_B <1$, there
exist positive parameters $\alpha,\beta, \kappa$ and $0<\rho <1$ depending on constants $\alpha_1,\alpha_2, \varLambda_1, \varLambda_2, \varLambda_3$ from 
Lemmas \ref{lm2}, \ref{lm3} and \ref{lm6} such that on any level $\ell\geq 0$,  the weighted term $\xi_\ell^2$ satisfies the following contraction property:
\begin{eqnarray}
 \xi_{\ell+1}^2 \leq \rho ~\xi_\ell^2, \label{contra}~~
 where ~~\xi_\ell^2:= \eta_\ell^2+\alpha e_{\ell}^2 +\beta \mu_\ell^2,\nonumber
\end{eqnarray}
whenever the initial mesh is chosen with small mesh-size. 
\end{theorem}\\
{\it Proof}.  
  For the {\bf Case (A)},  $\mu_\ell^2 \leq \kappa \eta_\ell^2$. The combination of (\ref{et_rd_11}) and (\ref{et_rd2})
 with a positive parameter $\beta$, to be chosen later, yields 
\begin{eqnarray}\label{con_2}
\eta_{\ell+1}^2+\beta \mu_{\ell+1}^2 \leq  \rho_A \eta_\ell^2  + (1+\delta_1) \beta \mu_\ell^2 +(\varLambda_2+\beta \varLambda_1 h_\ell^2) E_\ell^2,
\end{eqnarray}
where $h_\ell$ denotes the mesh-size at  the level $\ell$ of the triangulation.
With a choice of the initial mesh-size $h_3\leq\min \{ h_1,h_2, \frac{1}{\sqrt{\beta}}\}$,  multiply (\ref{er_re})
with the constant  $C_5=({\varLambda_1+ \varLambda_2 })/{\alpha_3}$, and then add 
with (\ref{con_2}) to obtain  
\begin{eqnarray*}
\eta_{\ell+1}^2+C_5(1-\alpha_1)e_{{\ell+1}}^2+\beta \mu_{\ell+1}^2 \leq  C_5 e_{{\ell}}^2+ \rho_A \eta_{\ell}^2  +( (1+\delta_1) \beta +C_5 \varLambda_3) \mu_\ell^2.
\end{eqnarray*}
Define $\alpha:=C_5(1-\alpha_1)$; $0<\alpha_4:=\frac{1}{2}\min\{1, \frac{1-\rho_A}{C_5C_{rel}} \}$,  
 and use the reliability result (\ref{estimator1}) that is, $e_{\ell}^2\leq C_{rel}(\eta_\ell^2+\mu_\ell^2)$ to obtain
\begin{eqnarray}\label{con_4}
\xi_{\ell+1}^2&&\leq C_5(1-\alpha_4)e_{\ell}^2+C_5\alpha_4 C_{rel}(\eta_\ell^2+\mu_\ell^2) +\rho_A \eta_\ell^2  +( (1+\delta_1) \beta +C_5 \varLambda_3) \mu_\ell^2 \nonumber \\
&&\leq C_5(1-\alpha_4)e_{\ell}^2+(\rho_A+C_5\alpha_4 C_{rel})\eta_\ell^2 +( (1+\delta_1) \beta +C_5 \varLambda_3+C_5\alpha_4 C_{rel}) \mu_\ell^2.
\end{eqnarray}
Since the marking criteria  implies  $E  \kappa \eta_\ell^2-E \mu_\ell^2>0$, where $E=2(C_5 \varLambda_3+C_5\alpha_4 C_{rel})$, an addition of this term on the right-hand side of (\ref{con_4}) yields
\begin{eqnarray}
\xi_{\ell+1}^2&&\leq C_5(1-\alpha_4)e_{\ell}^2+(\rho_A+C_5\alpha_4C_{rel}+E \kappa  )\eta_\ell^2\nonumber \\
&&~~~ +((1+\delta_1)\beta-(C_5 \varLambda_3+C_5\alpha_4 C_{rel}) ) \mu_\ell^2.
\end{eqnarray}
A use of Remark \ref{rmk2} and the definition   of $\alpha_4$ with the choice of the parameters   
$$\kappa <\kappa_0:= \frac{({1-\rho_A-C_5\alpha_4C_{rel}})}{E}, ~~\beta\geq 2(C_5 \varLambda_4+C_5\alpha_4C_{rel}), %
 ~~{\rm and}~~~ \delta_1<\min \Big \{ \frac{E}{2\beta}, \frac{\theta_B}{4-\theta_B} \Big { \} },$$
   yield on any level $\T_{\ell+1}$, 
 a contraction   $\xi_{\ell+1}^2 \leq  \rho_1 \xi_\ell^2 $, where 
\begin{eqnarray*}
0< \rho_1=\max \Big{\{ } \frac{C_5(1-\alpha_4)}{\alpha}, \rho_A+C_5\alpha_4C_{rel}+E \kappa , \frac{((1+\delta_1)\beta-E/2)}{\beta} \Big \}<1.
\end{eqnarray*}
%
For the {\bf Case (B)}: $\mu_\ell^2 > \kappa \eta_\ell^2$. 
Similar to proof of the {\bf Case (A)},  the  equation  corresponding to (\ref{con_4}) for the {\bf Case (B)} is 
\begin{eqnarray}\label{gg}
\xi_{\ell+1}^2&&\leq C_5(1-\alpha_4)e_{\ell}^2+C_5\alpha_4C_{rel}(\eta_\ell^2+\mu_\ell^2) +(1+\delta_2) \eta_\ell^2  +(\rho_B \beta +C_5 \varLambda_3) \mu_\ell^2\nonumber \\
&&\leq C_5(1-\alpha_4)e_{\ell}^2+(1+\delta_2+C_5\alpha_4C_{rel})\eta_\ell^2   +(\rho_B \beta +C_5 \varLambda_3+C_5\alpha_4C_{rel}) \mu_\ell^2.
\end{eqnarray}
The  marking  criteria  implies $  D \mu_\ell^2-D \kappa \eta_\ell^2>0$, where $0<D=\frac{1.5+C_5\alpha_4C_{rel}}{\kappa}$.
Add this term on the right-hand side of (\ref{gg}) to obtain 
\begin{eqnarray*}
\xi_{\ell+1}^2\leq C_5(1-\alpha_4)e_{\ell}^2+(2+C_5\alpha_4C_{rel}-D \kappa )\eta_\ell^2
+(\rho_B\beta  +C_5 \varLambda_3+C_5\alpha_4 C_{rel}+D ) \mu_\ell^2.
\end{eqnarray*}
A use of Remark \ref{rmk2}, the parameters choice and  $\beta >\frac{C_5 \varLambda_4+C_5\alpha_4C_{rel}+D}{1-\rho_B}$ 
 yields that, on any level $\T_\ell$  the contraction   $\xi_{\ell+1}^2 \leq  \rho_2 \xi_\ell^2 $  holds true, where $\rho_2$ is defined by 
\begin{eqnarray*}
0< \rho_2=\max \Big\{  \frac{C_5(1-\alpha_4)}{\alpha}, 2+C_5\alpha_4C_{rel}-D\kappa, \frac{\rho_B\beta+C_5 \varLambda_3+C_5\alpha_4C_{rel}+D }{\beta} \Big \}<1.
\end{eqnarray*}
 Finally, the combination of  both cases with 
    $0<\beta =2 \max\big{\{}C_5 \varLambda_3+C_5 \alpha_4 C_{rel},$ $ \frac{C_5\varLambda_3+C_5\alpha_4C_{rel}+D}{1-\rho_B}\big{\}},   $  $\rho=\max\{\rho_1,\rho_2\}$
     and the initial mesh-size $h_3$ implies that (\ref{contra}) holds. \hfill${\Box}$ 
\begin{theorem}
 (Convergence) Under the assumptions of Theorem \ref{thm1}, there exist a constant $\rho \in (0, 1)$ and $C_0>0$ 
   depending only on the given data and the initial triangulation such that 
\begin{equation*}
 \eta_\ell^2 +\alpha e_{\ell}^2 +\beta \mu_\ell^2\leq C_0 \rho^\ell.
\end{equation*}
\end{theorem}
{\it Proof}.  The proof is a consequence of the contraction property in Theorem \ref{thm1}.\hfill${\Box}$
\begin{remark}\label{rmk4} Using the reliability  result (\ref{estimator1}) and  the relation of  the error estimator $\eta_\ell^2$ and $\mu_\ell^2$ in the marking stratergy, the weighted term $\xi_\ell$, the  error  term $e_\ell$ and the error estimator terms $\eta_\ell$ and $\mu_\ell$ over triangulation $\mc T_\ell$ are obtained. 
Now $\xi_\ell$  satisfies
\[ \xi^2_\ell \leq \left\{
  \begin{array}{l l}
   C_6 \eta^2_\ell & \quad \text{for~~the~~  {\bf Case (A)}},\\
    \frac{C_6}{\kappa} \mu^2_\ell & \quad \text{for ~~the~~{\bf Case (B)}},
  \end{array} \right.\]
where $C_6=1+\alpha C_{rel}+ (\alpha C_{rel}+\beta)\kappa$. 
The reliability (\ref{estimator1}), Lemma \ref{lm8} with the  efficiency result (\ref{edeff}) implies
\begin{equation*}
\xi_\ell ^2\approx e_\ell^2 +\norm{h_{\mathcal T} f}^2.\end{equation*}
\end{remark}
\begin{lemma}\label{lm14} Let $\T_h$ and $\T_{H}$ be two nested triangulations. 
Then for small mesh-size $h_4>0$, it holds for  $0<h\leq h_4$
\begin{eqnarray}\label{1}
\xi_h^2\lesssim \xi_{H}^2. 
\end{eqnarray}
\end{lemma}
{\it Proof}.  
 From Lemma \ref{lm6}, Remark \ref{rmk3} and the efficiency result (\ref{lmeff}) implies
\begin{eqnarray}
\xi_h^2&&=\eta^2_h+\alpha e^2_h+\beta\mu_h^2\leq (C^{-1}_{\rm eff}+\alpha)e^2_{h} +2\beta \mu^2_{H} + \beta \varLambda_1 H ^2 E^2_{H}  \nonumber\\
&&\leq \frac{(C^{-1}_{\rm eff}+\alpha)}{1-\alpha_1}(e^2_{H}-(\frac{1}{2}-8C_3\varLambda_1 H^2) E^2_{H}+\varLambda_3 \mu^2_H )+2\beta\mu^2_{H} +\beta\varLambda_1 H^2 E_{H}^2.\label{2}
\end{eqnarray}
Select the initial mesh-size $h_4>0$ such that the coefficient  of $E_H^2$ is non-positive for $0<h\leq h_4.$
 Then (\ref{2}) implies (\ref{1}). Note that the inequality constant in (\ref{1}) is independent of AMFEM marking parameters. This concludes the proof. \hfill${\Box}$
\subsection{Quasi-optimality}
In this subsection, the quasi-optimal convergence \cite{steven} of the adaptive algorithm MFEM is discussed with the help of the quasi-discrete reliability and the contraction property. 
\begin{definition} \cite{Car_apo, steven}(Approximation class) Given an initial triangulation $\mathcal T_0$ of $\Omega$ and $s>0$, the approximation class is defined as
\begin{eqnarray*}
\mc A_s :=\{ ({\bf p},u,f)\in H({\rm div},~ \Omega) \times L^2(\Omega) \times L^2(\Omega)~|~ \norm{({\bf p},u,f)}_{\mc A_s} <\infty \} ~~~{\rm with}
\end{eqnarray*}
\begin{eqnarray*}
 && \norm{({\bf p},u,f)}_{\mc A_s} := \sup_{N\in \mathbb N}\Big( N^s \inf_{|\mc T|-|\mc T_0|\leq N} (e^2(\mc T)+\norm{h_{\mathcal T} f}^2)^{{1}/{2}}\Big) {\rm and} ~ e^2(\T)=e^2_p(\T)+e^2_u(\T).
\end{eqnarray*}
Here, the infimum is  over all regular and NVB-generated refinements $\mc T$ of $\mc T_0$
 ~with the number of element domains $|\mc T|\leq N+|\mc T_0|$ and the exact error $e(\mc T)$ in MFEM solution.
\end{definition}
%
An adaptive mixed finite element method is quasi-optimal convergent in the sense that 
given $({\bf p}, u, f)\in \mathcal A_s,~ {\rm and} ~\J \in \mathbb N$, the  AMFEM  algorithm  generates a triangulation $\mathcal T_\J$ with discrete  solution $({\bf p_\J},u_\J)\in RT_0(\mathcal T _\J)\times P_0(\mathcal T_\J)$ such that 
$$|\mathcal T_\J|-|\mathcal T_0| \leq \xi_\J ^{-\frac{1}{s}}\approx (e_\J^2 +\norm{h_{\mathcal T}f}^2)^{-\frac{1}{2s}}.$$
\begin{theorem}(Quasi-optimality) Assume $({\bf p}, u, f)\in \mathcal A_s.$ Let $\{\mc T_\J \}_{\J\geq 0}$ be the sequence of the meshes generated  by AMFEM algorithm and $\{ ({\bf p_\J},u_\J)\in {RT_0(\mc T_\J)}\times P_0(\mc T_\J)\}_{\J\geq 0}$ be a corresponding sequence of approximate solutions.
Then, for   small initial mesh-size $h_0$, the following estimates holds true for $0<h\leq h_0$
\begin{equation}\label{obtim}
|\mathcal T_\J|-|\mathcal T_0| \lesssim \xi_\J^{-\frac{1}{s}}. 
\end{equation}
\end{theorem}
{\em Proof}.
Consider
 \begin{eqnarray}
\qquad \qquad |\mathcal T_\J|-|\mc T_0| \leq \sum_{\ell=0}^{\J-1} \Big(|\mc T_{\ell+1}|-|\mc T_\ell| \Big) 
\end{eqnarray}
Now $|\T_{\ell+1}|-|\mc T_\ell|\lesssim |\mc M_\ell|$ for $\ell\geq 0$ \cite{biven,steven},
~ where $\T_{\ell+1}$ is a refinement of $\T_\ell$ and $\mc M_\ell$ denotes the set of the marked edges or elements in the triangulation level $\T_\ell$. Then
\begin{eqnarray}\label{qo1}
|\mathcal T_\J|-|\mc T_0| 
&\lesssim& \sum_{\ell=0}^{\J-1} |\mc M_\ell|. 
\end{eqnarray}
To estimate of $|\mc M_\ell|$, use the characterization of the approximate class and the overlay.\\
If $({\bf p}, u, f)\in \mathcal A_s$, then for $\e_1:=\tau_1 \xi_\ell$, 
there exists some admissible triangulation $\mc T_{\e_1}$ obtained as refinement of the initial triangulation $\T_0$  such that 
\begin{eqnarray}\label{opt1}
\xi^2_{\e_1}=\eta^2_{\e_1}+\alpha e^2_{\e_1}+\beta \mu^2_{\e_1}\leq \e_1^2 ~~~\rm{and}~~~ |\T_{\e_1}|-|\T_0|\lesssim \e_1^{-\frac{1}{s}}.
\end{eqnarray}
Here, $0<\tau_1=\min\{ \tau_2,\tau_3 \}$, and $\tau_2 ~{\rm and} ~\tau_3$ will be specified later.\\
Let $\mc T_{\ell+{\e_1}}:=\mc T_{\e_1}\oplus \T_\ell $ be the overlay of $\T_{\e_1}$ and $\T_\ell$. As in  \cite{quasi}, the  number of elements of the overlay $\mc T_{\ell+{\e_1}}$  can be bounded by
\begin{eqnarray}\label{opt2}
|\mc T_{\ell+{\e_1}}|-|\T_\ell|\leq|\T_{\e_1}|-|\T_0|.
\end{eqnarray}
To estimate  $|\mc M_\ell|$ with the first case of Mark algorithm, that is, $\kappa \eta^2_\ell \geq \mu^2_\ell$,
first define $\mc M^\dagger :=\mc E_\ell\setminus \mc E_{\ell+\e_1} \subset \mc E_\ell$ as the set of edges of $\T_\ell$ being refined in $\T_{\ell+\e_1}$. Note that, if  $\mc M^\dagger$ satisfies the marking criteria of the {\bf Case (A)}, that is, 
\begin{eqnarray}\label{opti_2}
\theta_A \eta_\ell^2 \leq \eta^2_\ell(\mc M^\dagger),
\end{eqnarray} 
then,  $|\mc M_\ell|\leq|\mc M^\dagger|$, where $|\mc M_\ell|$ be the set of marked edges at level $\ell$.\\
The quasi-discrete reliability result Lemma \ref{coro2} over the triangulations $\T_\ell$ and $\T_{\ell+\e_1}$  with  $\epsilon^2=1/4C_4 $ and the initial mesh-size $h_5$ shows
\begin{eqnarray*}
C_4 \eta_\ell^2 (\mc M^\dagger)\geq  E_\ell^2 -C_4 osc^2_\ell (f)-\frac{1}{4}(e^2_\ell+e^2_{{\ell+\e_1}})-C_4(\mu^2_\ell+\mu^2_{\ell+\e_1}).
\end{eqnarray*}
Remark \ref{rmk3} yields $\mu^2_{{\ell+\e_1}} \leq 2 \mu_\ell^2+\Lambda_1 h_\ell^2 E^2_\ell$, where $h_\ell$ denotes the mesh-size over the  triangulation $\T_\ell$ and note that $ osc^2_\ell (f)\leq \mu_\ell^2$. Hence,  
\begin{eqnarray}\label{qo3}
C_4 \eta_\ell^2 (\mc M^\dagger)\geq(1-C_4\varLambda_1 h_\ell^2)  E_\ell^2 -4 C_4 \mu_\ell^2 -\frac{1}{4}(e^2_\ell+e^2_{{\ell+\e_1}}).
\end{eqnarray}
For  small initial mesh-size $h_4$, Lemma \ref{lm14}, (\ref{opt1}), the choice of $\epsilon_1$ and Remark \ref{rmk4} imply
 $$e_{{\ell+\e_1}}^2 \leq \xi_{\ell+\e_1}^2\lesssim  \xi^2_{\e_1}\lesssim  \tau_1^2 \xi^2_\ell \lesssim \tau_1^2  \eta_\ell^2 \leq  \tau_2^2 \eta_\ell^2,$$ 
 that is, $e_{{\ell+\e_1}}^2 \leq C_7 \tau_2^2 \eta_\ell^2$.  Thus, Corollary \ref{coro3} supplies the lower bound for $E_\ell$. With these results, (\ref{qo3}) leads to
\begin{eqnarray}
C_4 \eta_\ell^2 (\mc M^\dagger)&\geq &\frac{1}{2} (\frac{1}{2}-C_4\varLambda_1^2 h_\ell^2) e_\ell^2-\frac{5}{4}C_7 \tau^2 \eta_\ell^2-\big( 4C_4+ \frac{\varLambda_3}{2}\big)\mu_\ell^2.
\end{eqnarray}
Note that some positive terms are neglected from the right-hand side. 
Choose the initial mesh-size $h_6 :=\min\{h_2,h_5,\frac{1}{(4C_4\varLambda_1)^{1/2}} \}$ and  the marking parameter $\theta_A=\min \{1 , \frac{c_{\rm eff}}{9C_4}\}$.  
 The {\bf Case (A)} relation $\kappa \eta_\ell^2\geq \mu_\ell^2$ along-with  the efficiency result (\ref{edeff}) leads to
\begin{eqnarray}
C_4 \eta_\ell^2 (\mc M^\dagger)&\geq &\Big(\frac{ C_{\rm eff}}{8} -\frac{5}{4} C_7\tau_2^2 - \kappa\big( 4C_4+ \frac{\varLambda_3}{2}\big)-C_4 \theta_A\Big)\eta_\ell^2+ C_4 \theta_A \eta_\ell^2. \nonumber
\end{eqnarray}
The  selection of $ \kappa < \min \{ \kappa_0,  \frac{{C_{\rm eff}}/{8}-C_4\theta_A }{4C_4+ {\varLambda_3}/{2}}\big\}$, and $\tau_2^2:=\big( \frac{(c_{\rm eff}/8)-\kappa ( 4C_4+ {\varLambda_3}/{2})-C_4 \theta_A}{(5/4)C_7} \big) $ imply that (\ref{opti_2}) holds.\\
Since  $\mc M_\ell$ is chosen to be  the minimal cardinality set satisfying (\ref{opti_2}),   Lemma 4.4 of \cite{Car_apo},  
(\ref{opt2}) and (\ref{opt1}) altogether imply
\begin{eqnarray}\label{opti_3}
|\mc M_\ell|\leq |\mc M^\dagger_\ell| \lesssim |\T_{\ell+\e_1}| -|\mc T_\ell |\leq|\T_{\e_1}|-|\T_0|\lesssim\e_1^{-\frac{1}{s}}\lesssim\xi^{-\frac{1}{s}}_\ell.
\end{eqnarray}
Consider the {\bf  Case (B)} of Mark algorithm, that is, $\kappa \eta^2_\ell \leq \mu^2_\ell $.  Let $\mc M^* = \T_\ell\setminus \T_{\ell+\e_1}\subset \T_\ell$ be the set of elements of $\T_\ell$ refined in $\T_{\ell+\e_1}$. The proof of $\mc M^*$ satisfies the marking criteria of the {\bf Case (B)}, that is,
\begin{eqnarray}\label{opti_4}
\theta_B \eta^2_\ell\leq \eta^2_\ell(\mc M^*),
\end{eqnarray} 
and this will imply $|\mc M_\ell|\leq|\mc M^*|$, where $\mc M_\ell$ is the set of the marked 
elements at the level $\ell$.\\
Corollary \ref{lm13} and   the quasi-discrete reliability result in Lemma \ref{coro2}, with  $\epsilon=1 $ and initial mesh-size $h_7$ over the nested triangulations $\T_{\ell+\e_1}$ and $\T_\ell$ imply 
\begin{align*}
\mu^2_{{\ell+\e_1}} &\geq \frac{\gamma_0^2}{2} \mu^2_\ell -2\norm{h_\ell(f_\ell -f_{{\ell+\e_1}})}^2-C_{1} h_\ell^2 E_\ell^2 \nonumber\\
 &\geq \big(\frac{\gamma_0^2}{2}-C_4C_{1} h_\ell^2\big) \mu^2_\ell -(2+C_4C_{1}h_\ell^2)\norm{h_\ell(f_\ell -f_{{\ell+\e_1}})}^2\nonumber\\
&~~~ -C_{1}C_4 h_\ell^2\big(\eta^2_\ell +{}(e_\ell^2+e^2_{{\ell+\e_1}})
+\mu^2_{\ell+\e_1}\big).
\end{align*}
The reliability result, the  relation $\kappa \eta_\ell^2 <\mu_\ell^2$ for the  {\bf Case (B)} and the rearrangement  of  terms   imply 
\begin{eqnarray}\label{op2}
~~\mu^2_{{\ell+\e_1}}\geq\frac{(\frac{\gamma_0^2}{2}-C_{8} h_\ell^2 )}{(2+C_4C_{1} h_\ell^2)}\mu^2_\ell -\norm{h_\ell(f_\ell -f_{{\ell+\e_1}})}^2- {C_{1}C_4 h_\ell^2}e^2_{{\ell+\e_1}},
\end{eqnarray}
where $C_{8}:= C_4C_{1}(Crel+1)(1+1/\kappa)$.
 Lemma \ref{lm14}, the choice of $\epsilon_1$, that is, $\epsilon_1=\tau_1\xi_\ell$ and   the  marking criteria    $\kappa \eta_\ell^2< \mu_\ell^2$ in Remark {\ref{rmk4}} result in
\begin{eqnarray}\label{op1}
\alpha e_{{\ell+\e_1}}^2+\beta \mu^2_{{\ell+\e_1}}\leq\xi^2_{{\ell+\e_1}} \lesssim \xi_{\e_1}^2\lesssim \tau_1^2 \xi^2_\ell \lesssim \tau_3^2 \mu_\ell^2,
\end{eqnarray}
that is, $\alpha e_{{\ell+\e_1}}^2+\beta \mu^2_{{\ell+\e_1}}\leq C_9\tau_3^2 \mu_\ell^2 $. The combination of (\ref{op2})-(\ref{op1}), for  small  initial mesh-size 
$h_8=\min\big{\{ }h_7, h_4, ({\alpha}/{\beta C_{1}C_4})^{1/2},({\gamma_0^2}/{4C_{8})})^{1/2}\big{\}}$ and some simplifications show
\begin{equation*}
\beta{\mu^2(\mc M^*)}\geq \beta\norm{h_\ell(f_\ell -f_{{\ell+\e_1}})}^2\geq \Big[\frac{\beta{\gamma_0^2}}{4(2+C_4C_{1} h_\ell^2)}-2C_9 \tau_3^2  
- \beta \theta_B\Big] \mu^2_\ell+ \beta \theta_B \mu^2_\ell.
\end{equation*}
The selections $ \tau^2_3: =\frac{\beta{\gamma_0^2}}{16(2+C_4C_{1} h_\ell^2)C_9}$ and $ \theta_B =\frac{{\gamma_0^2}}{8(2+C_4C_{1} h_\ell^2)}$ lead to (\ref{opti_4}).
Since $\mc M_\ell$ is chosen to be the minimal cardinality set satisfying (\ref{opti_4}),  (\ref{opt2}) and (\ref{opt1}) yield
\begin{eqnarray}\label{opti_5}
|\mc M_\ell|\leq|\mc M^*|\leq |\T_{\ell+\e_1}|-|\T_\ell|\leq|\T_{\e_1}|-|\T_0|\lesssim\e_1^{-\frac{1}{s}}\lesssim\xi^{-\frac{1}{s}}_\ell.
\end{eqnarray}
Now a combination of both the cases, that is, (\ref{opti_3}) and (\ref{opti_5}) with (\ref{qo1}) leads to 
\begin{eqnarray}
|\mathcal T_\J|-|\mc T_0| 
&\lesssim& \sum_{\ell=0}^{\J-1} |\mc M_\ell| 
\lesssim \sum_{\ell=0}^{\J-1} \xi_\ell^{-\frac{1}{s}}.\nonumber 
\end{eqnarray}
A use of  the contraction property in Theorem \ref{thm1} shows
\begin{align}
|\mathcal T_\J|-|\mc T_0| 
&\leq ( \rho^{\J(\frac{1}{2s})}+ \rho^{(\J-1)(\frac{1}{2s})}+\cdots+\rho^{\frac{1}{2s}} \big) \xi^{-\frac{1}{s}}_{\J}\nonumber\\
 &\leq  \xi^{-\frac{1}{s}}_{\J} \sum_{j=1}^{\J} \rho^{\frac{j}{2s}}
\lesssim \xi_{\J}^{-\frac{1}{s}}\Big[ \frac{1}{1-\rho^\frac{1}{2s}}\Big] \lesssim  \xi^{-\frac{1}{s}}_{\J},
\end{align}
 and this concludes the proof. \hfill${\Box}$
\section{Numerical Experiments}
This section shows the performance of the adaptive algorithm (AMFEM) on some benchmark problems. 
\begin{example} \label{example1} For the problem (\ref{eq1}), set $\Omega=(0,1)\times (0, 1)$ and the coefficients ${\bf A}=I,{\bf b}=(1,1),\gamma=2$. Choose the right-hand side $f$ and $u|_{\partial \Omega}$  such that $u=exp(-100\norm{x-x_0}^2)$.
%
\end{example}\\
 The initial uniform criss-cross triangulation  $\T_0$ has mesh-size $h=0.25$ as shown in Fig. 5.1(a). The adaptive algorithm  is performed  with parameters  $\theta_A=0.5, \theta_B=0.5$ and  $\kappa=0.8$.
\begin{figure}[!ht]
\centering
\subfloat[] 
{\includegraphics[width=6.0cm]{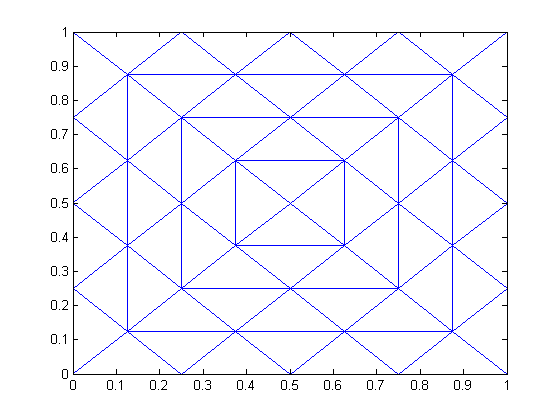}} \subfloat[ ]{\includegraphics[width=6.0cm]{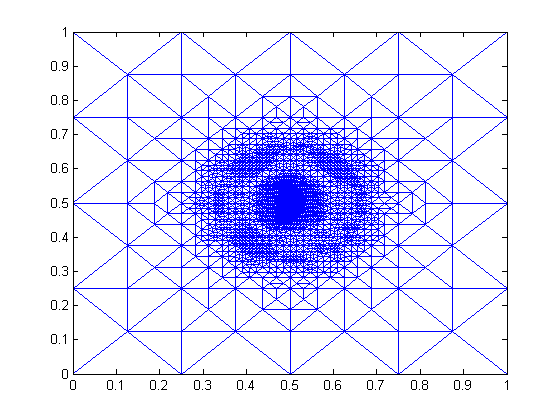}}
\caption{(a) The initial triangulation (b) The adaptive refined mesh at the level 10.}
\end{figure}
\begin{table}[h!]
\footnotesize
\begin{center}
  \begin{tabular}{ c | c | c  |c|c |c}
    \hline
   Ndof &$\norm{u-u_\ell}$ &$\norm{{p-p}_\ell}$&$\eta_\ell$&$\mu_\ell$& ~Marking Case\\
 \hline \hline
 168    &  0.0571&    1.3649 &   1.5499 &   4.8186  & B\\
 203   &   0.0435 &   0.9670  &  2.7266  &  4.1602  &  B\\
 233    &  0.0336  &  0.8209   & 2.9629   & 3.3392   & B\\
368   & 0.0254&    0.6160&    2.2798&    2.3559&   B\\
523    &0.0197&    0.4471 &   1.5828 &   1.6907 &  B\\ \hline
948   & 0.0152 &   0.3287  &  1.1607  &  1.1389  &  B\\
1783   & 0.0103 &   0.2451  &  0.8943  &  0.8198  &  B\\
3373    &0.0074  &  0.1647   & 0.6013   & 0.5797   & B\\
 6138   &0.0058   & 0.1334&    0.4908&    0.4182&    A\\
6918    &0.0048    &0.1027 &   0.3692 &   0.3997 &   B\\ \hline
12053   & 0.0041&    0.0866 &   0.3120 &   0.2941 &   B\\
21298  &  0.0033 &   0.0698  &  0.2527  &  0.2177  &  A\\
 23198&   0.0026  &  0.0547   & 0.1991   & 0.2137   & B\\
 40658   &0.0022   & 0.0461&    0.1685&    0.1599&    B\\
  \hline
  \end{tabular}
\end{center}
    \caption{Numerical results of AMFEM for Example \ref{example1}}
\end{table}
\begin{figure}[!ht]
\centering
{\includegraphics[width=6.2cm]{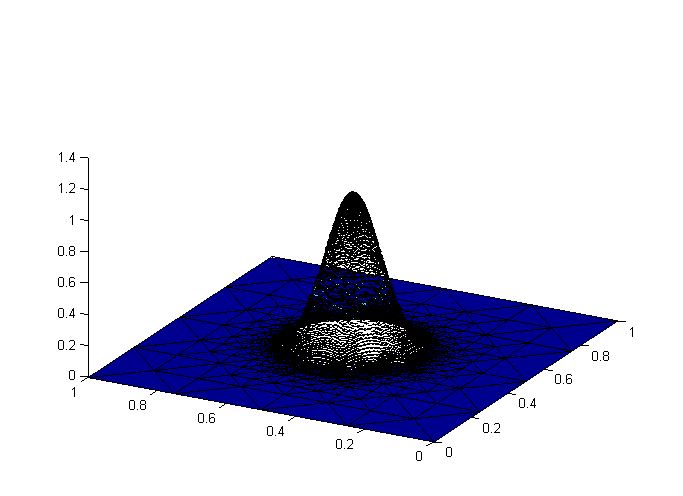}}
{\includegraphics[width=6.2cm]{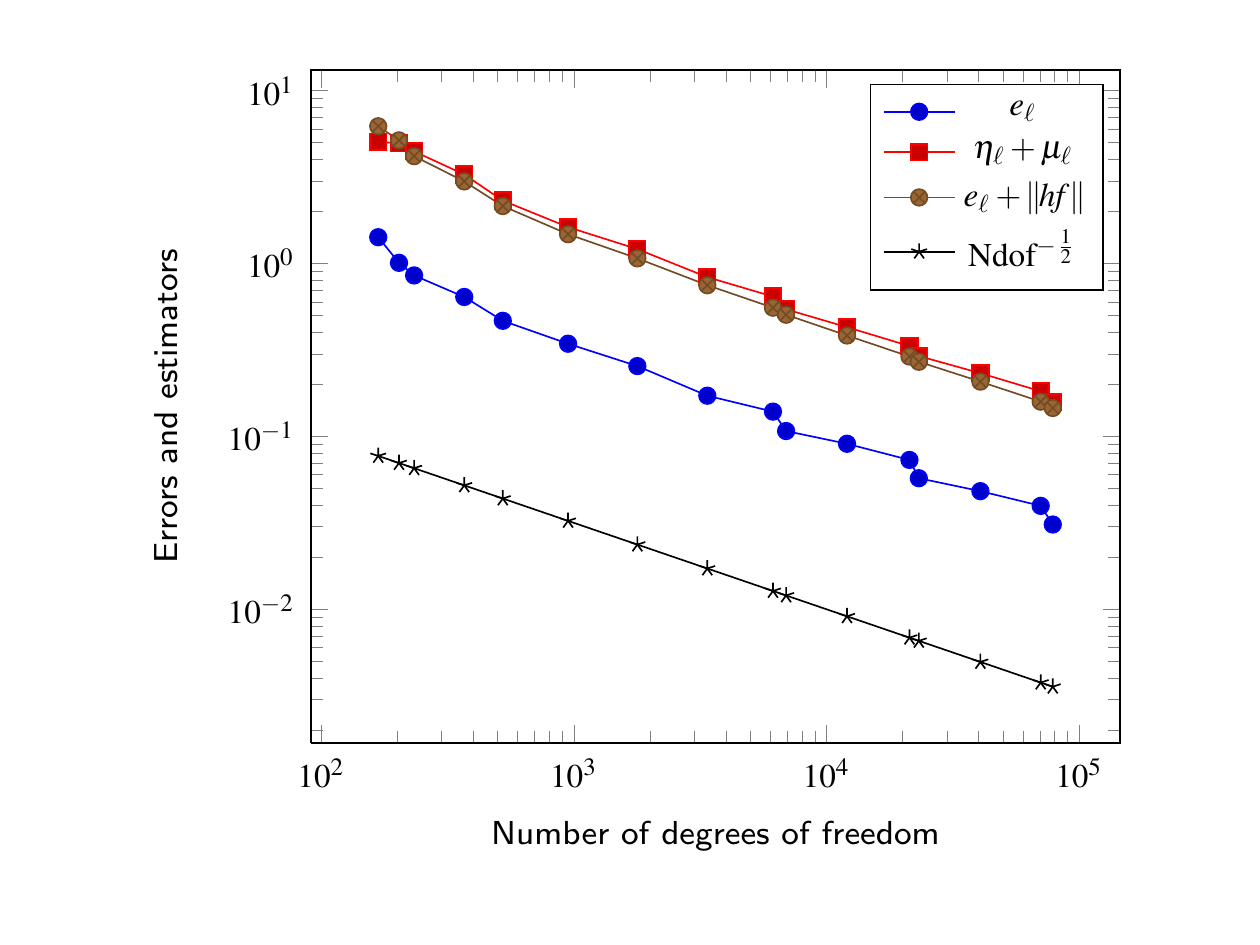}}

\caption{Discrete solution and convergence plot}
\end{figure}
Table 5.1 displays the experimental results for  errors $\norm{u-u_\ell}$ and  $\norm{{\bf p-p}_\ell}$, the edge and volume estimators $\eta_\ell$ and $\mu_\ell$,  from (\ref{eta}) and (\ref{mu})
for several consecutive levels $\ell$  of AMFEM with the number of degrees of freedom and the marking {\bf Case (A) or (B)}. 
For this example, the load function $f$ exhibits a relatively large
variation in the domain and hence, within the elements. 
It is observed that  the volume estimator $\mu_\ell$ is more than the edge-estimator $\eta_\ell$ in  several levels due to large  data oscillations, which results in the use of the  marking {\bf Case (B)} for most of the refinement levels. The reduction of $\mu_\ell$ influences the reduction in the edge-estimator and both lead to the optimal convergence rate ${\rm Ndof}^{-1/2}$ for errors, where Ndof denotes the number of degrees of freedom.
Fig. 5.1(b) shows the adaptive refined mesh by the AMFEM algorithm at the level 10, where the number of degrees of freedom is 6918. Fig. 5.2 depicts the approximate solution $u_h$ and summarises the convergence rates for variables ${\bf p}, u$, and sum of the estimators with respect to 
the number of degrees of freedom. 
\begin{example}\label{example2}
For the {\it Crack problem} \cite{3matlab}, consider the PDE (\ref{eq1}) with coefficients $\mathbf A=I, ~{\bf b}=(x-1,y+1)$ and $\gamma=4$ on
$~\Omega =\{(x,y)\in \mathbb{R}^2:|{{\bf x}}|\leq 1\setminus[0, ~1]\times \{0\}\}$ with Dirichlet boundary condition and
  exact solution $~u(r,\theta)=r^{{1}/{2}} \sin{\theta}/{2} -{r^2}/{2} \sin ^2(\theta)$. 
\end{example}   
The initial uniform  triangulation  $\T_0$ has  mesh-size $h=0.25$. The adaptive algorithm  is performed with parameters  $\theta_A=0.3, \theta_B=0.3$ and  $\kappa=1$.    
\begin{table}
\footnotesize
\begin{center}
  \begin{tabular}{ c | c | c  |c|c |c}
     \hline
   Ndof &$\norm{u-u_\ell}$ &$\norm{{p-p}_\ell}$&$\eta_\ell$&$\mu_\ell$& Marking Case\\
 \hline \hline
46 &0.0880&    0.4021 &   2.0875 &   2.5734   & A\\
 56  & 0.0899 &   0.3979 &   3.4950 &   1.9191   & B\\
98 & 0.0771  &  0.3557  &  2.3290   & 1.8010  &  B\\
161 &   0.0633   & 0.3085    &1.3106    &1.4866 &   A\\
207  &  0.0556    &0.2643   & 1.1180    &1.0874  & B\\ 
383   & 0.0429  &  0.2401   & 0.6512    &0.8751   & A\\ \hline
469    &0.0371   & 0.2042  &  0.5813    &0.6579    &A\\
628  &  0.0308   & 0.1695 &   0.4662    &0.4896    &A\\
901  &  0.0269&    0.1400&    0.3579 &   0.3431 &   B\\
1436  &  0.0207&    0.1225&    0.1971 &   0.2704 &   A\\  
1821   & 0.0183 &   0.1023 &   0.1762  &  0.1949  &  A\\
2618   & 0.0150  &  0.0844  &  0.1433   & 0.1340   & B\\ \hline
3957   & 0.0123   & 0.0746   & 0.0795    &0.1086    &A\\
5162    &0.0111&    0.0603    &0.0669&    0.0698    &A\\
7482    &0.0092 &   0.0486&    0.0528 &   0.0453    &B\\  
12145    &0.0074 &   0.0434&    0.0282 &   0.0384&    A\\
15726&    0.0065  &  0.0351 &   0.0231  &  0.0246 &   A\\
22520 &   0.0054   & 0.0280  &  0.0169   & 0.0158  & B\\ \hline 
           \hline
  \end{tabular}
\end{center}
    \caption{Numerical results of AMFEM for Example \ref{example2}}
\end{table}

\begin{figure}[!ht]
\centering
{\includegraphics[width=6.2cm]{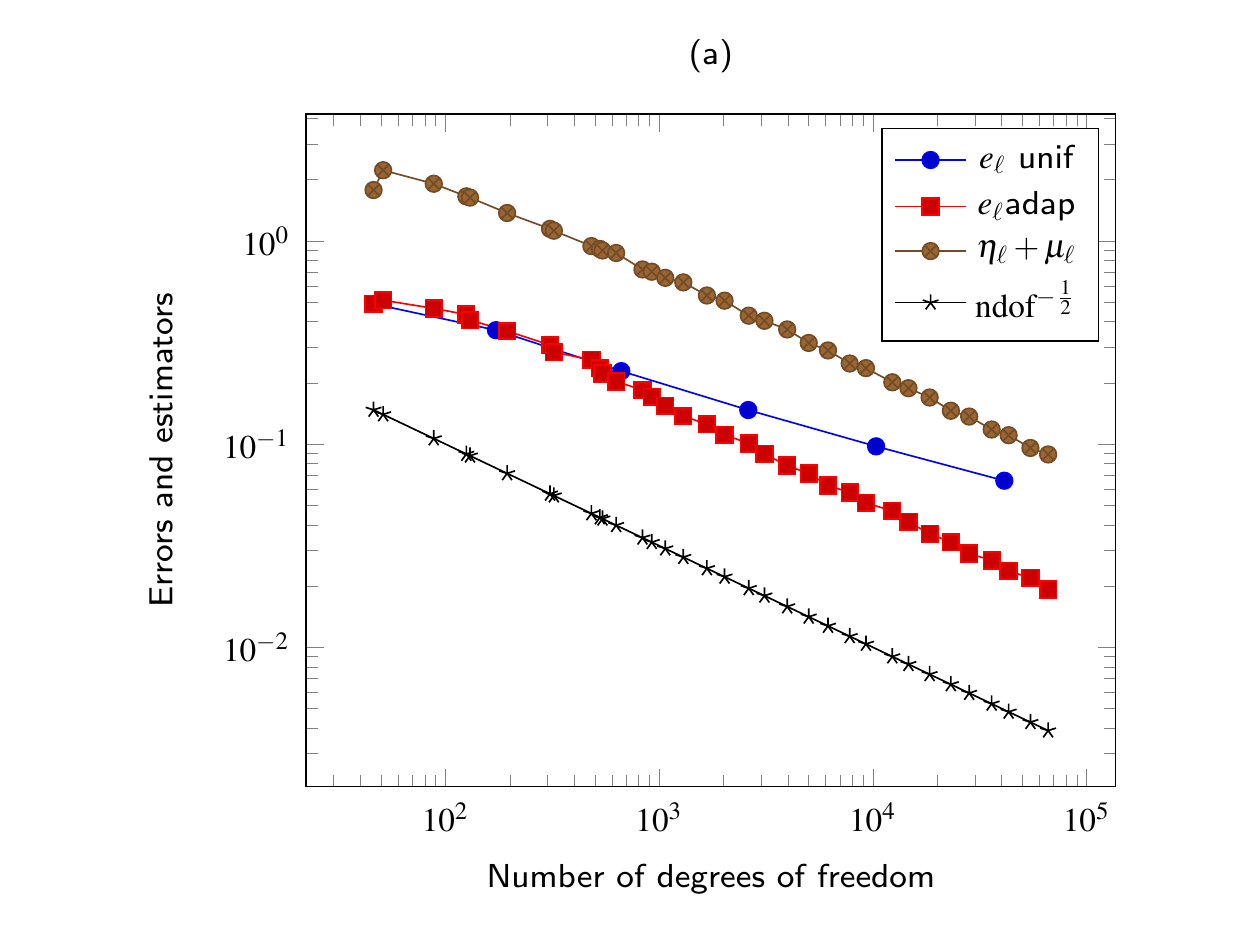}}
{\includegraphics[width=6.2cm]{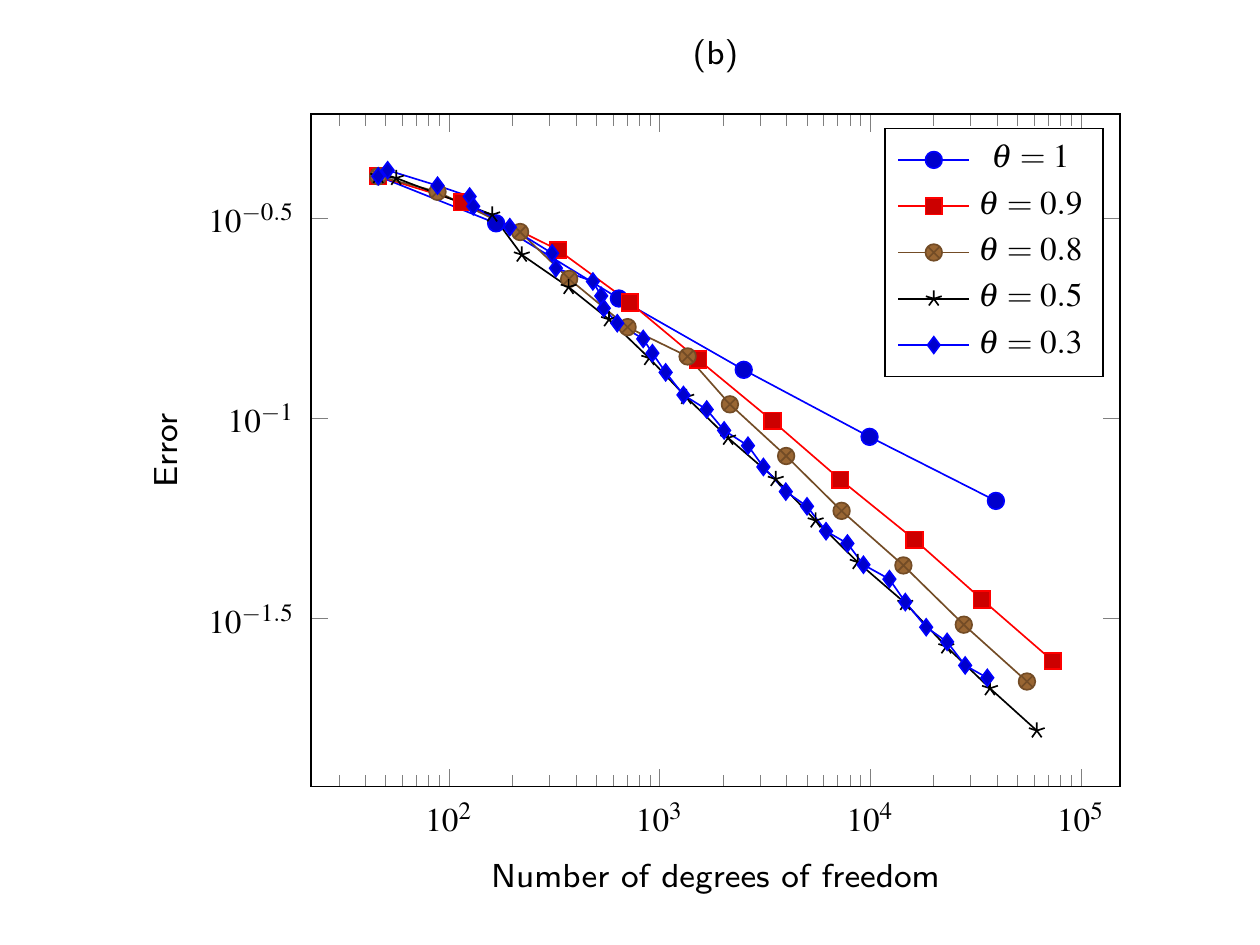}}

\caption{(a) Convergence plots for $\theta_A,\theta_B=0.3$ (b) Error convergence for $\theta=\theta_A=\theta_B$ equal to 0.3, 0.5, 0.8, 0.9, 1}
\end{figure}     
       Table 5.2 displays the experimental results for the errors  $\norm{u-u_\ell}$ and $\norm{{\bf p-p}_\ell}$, the estimators $\eta_\ell$ and $\mu_\ell$ 
for several consecutive levels $\ell$  of AMFEM with the number of degrees of freedom and the marking {\bf Case (A) or {\bf (B)} } for this problem. 
The right-hand function $f$ of (\ref{eq1}) in this example  is not smooth and the solution has singularity at the origin, and hence  the adaptive algorithm utilizes both the marking cases to achieve the optimal convergence rate. Figure 5.3(a) shows  optimal convergence rate ${\rm Ndof}^{-1/2}$ for adaptive algorithm while uniform refinement yields a suboptimal convergence  rate ${\rm Ndof}^{-1/4}$ for the above mentioned parameters. Figure 5.3(b) shows the convergence rate of the error  for  different values of $\theta_A$ and $ \theta_B$. This supports the theoretical prediction of the choice of marking parameters $\theta_A$ and $\theta_B$ in the proof of Theorem 4.5, that is, the quasi-optimality can be 
achieved only for  small values of  $\theta_A~ {\rm and}~ \theta_B$.
\begin{example}\label{example3} 
Consider the PDE (\ref{eq1}) with  coefficients  $~ \mathbf A=I,~~\mathbf {b} =(0, 0),$~ 
$\gamma=-8.9$,  Dirichlet boundary condition on the L-shaped domain
$\Omega =(-1,1)\times (-1,1)\setminus [0,1]\times [-1,0]$ and the exact solution  given in polar coordinates as $u(r,\theta)=r^{2/3} \sin\big( 2\theta/3 \big).$
 \end{example}    
\begin{figure}[!ht]
\centering
{\includegraphics[width=6.0cm]{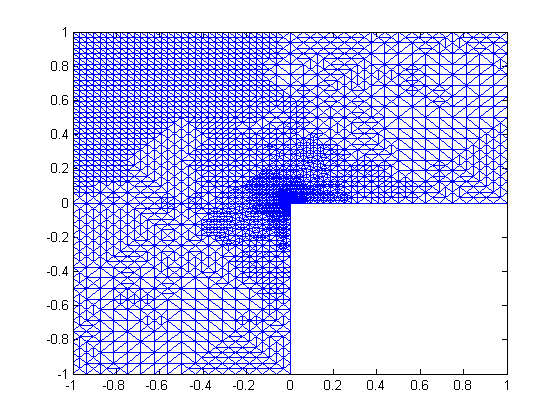}} 
{\includegraphics[width=6.2cm]{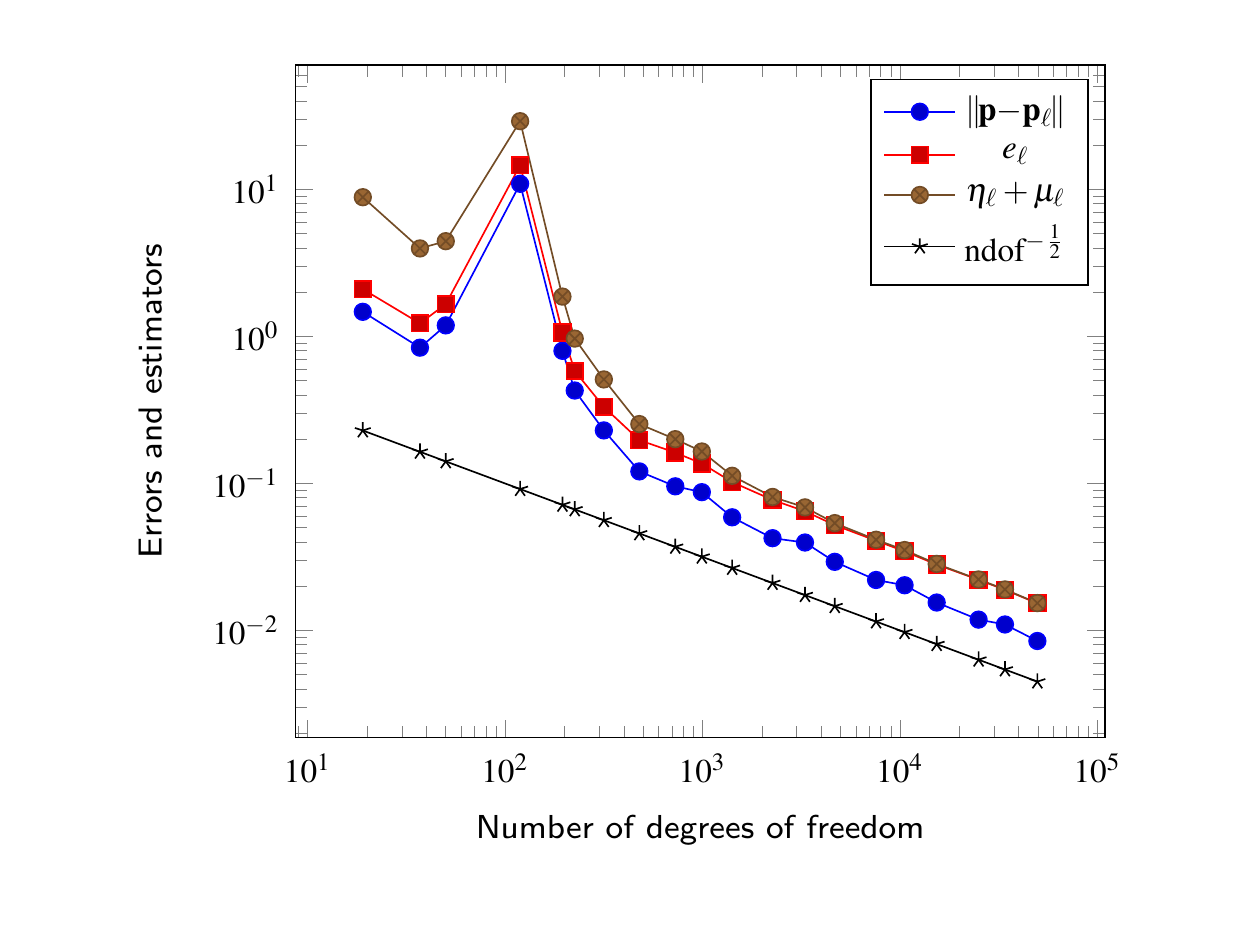}} 
%
\caption{Adaptive mesh-refinement and the convergence plot}
\end{figure}     
The numerical experiment  is performed with initial mesh-size $h=0.25$ and uniform triangulation. The adaptive algorithm  is performed with parameter choice $\theta_A=0.5, \theta_B=0.5$ and  $\kappa=2$.
Since  the  forcing function $f$ is smooth, the oscillation term has higher order convergence than the remaining terms in the estimator. In this example, the marking is based on the edge estimator  to obtain the optimal convergence. 
 The refinement due to the data oscillations and the volume estimators  can be avoided since they play a minor role in the convergence. 
        Figure 5.4 displays the adaptive mesh-refinement and the optimal convergence rate ${\rm Ndof}^{-1/2}$  of the errors and the estimators for a sufficiently small mesh-size $h$.
\subsection*{5.4.  Observations}
This subsection deals with a few observations.
\begin{itemize}
 \item If the given function $ f$ has large variation in the domain and within the elements, 
 then adaptive algorithm chooses 
    {\bf Case (B)} marking criteria and achieves the convergence (see 5.1).
\item When $f$ is not smooth and the solution also has singularity, then the adaptive algorithm utilizes both the marking cases to achieve the optimal convergence rate (see 5.2).
\item For the  smooth function $f$, the oscillation has higher order convergence than the remaining  estimator terms.
 In such cases, the data oscillations and the  volume estimators 
 have a very minor role in the convergence. Thus, the marking would be based on the edge-based error estimator to capture singularity of the solution (see 5.3).  
\end{itemize}
\section{Conclusions} In this work, the convergence and the quasi-optimality  of adaptive mixed finite element method
is analysed for the non-symmetric and indefinite second order elliptic equations using the lowest order Raviart-Thomas elements. The adaptive algorithm in Subsection 2.1 is designed as a combination of  the edge and the volume error estimators. The numerical experiments confirm the efficiency of this algorithm and  support the theoretical findings.
\subsection*{Acknowledgements} The first author acknowledges the financial support of  Council of Scientific and Industrial Research (CSIR), Government of India. The authors sincerely thank Professor Carsten Carstensen,  Humboldt
University, Berlin  for his constructive comments and suggestions.

{\footnotesize
  \begin{center}
  {\bf Appendix I}\\
{List of the constants}
\begin{longtable}{| c | p{6cm} |  p{2.8cm}|p{1.1cm}|}

\hline \multicolumn{1}{|c|}{\textbf{Constant}} & \multicolumn{1}{p{6cm}|}{\textbf{Dependency on }} & \multicolumn{1}{p{2.8cm}|}{\textbf{Appears first in}} & \multicolumn{1}{p{1.1cm}|}{\textbf{Value }} \\ \hline 
\endfirsthead

\multicolumn{4}{c}%
{{\bfseries \tablename\ \thetable{} -- continued from previous page}} \\
\hline \multicolumn{1}{|c|}{\textbf{Constant}} &
\multicolumn{1}{p{6cm}|}{\textbf{Dependency on }} &
\multicolumn{1}{p{2.8cm}|}{\textbf{Appears first in}}  &
\multicolumn{1}{p{1.1cm}|}{\textbf{Value }} \\ \hline 
\endhead

\hline \multicolumn{4}{|r|}{{Continued on next page}} \\ \hline
\endfoot

\hline \hline
\endlastfoot
&{\bf I. Natural Constants~~~~} && \\ \hline
$ C_{\rm rel}$   & Coefficients of (1.1) and interpolation  &    &  Positive \\
&constants&Theorem 3.2 &\\ \cline{1-2} \cline{4-4}
$\epsilon$   & Coefficients of (1.1) and interpolation  & ({\it A~posteriori }   &  Positive \\
&constants& {\it estimates}) &\\  \cline{1-2} \cline{4-4}
$h_1$& $\epsilon$& &  Positive \\ \hline
$  C_{\rm eff} $ & Coefficients of (1.1),  inverse inequality  &   Lemma 3.3 &  Positive \\  
&and finite overlap&({\it efficiency})  &\\ \cline{1-3} \cline{3-4}
$ C_{1}$   & Coefficients of (1.1)&   Lemma 3.4 ({\it Volume} &  Positive \\
&&{\it  estimator reduction})& \\ \cline{1-3}\cline{4-4}
$  C_{2} $ & Coefficients of (1.1),  inverse inequality  &   Lemma 3.6 ({\it edge}&  Positive \\ 
&and finite overlap&({\it  estimator reduction})& \\\cline{1-4}
$ C_{3}$   & Coefficients of (1.1)  &  Lemma 3.7 ({\it Quasi-} &  Positive \\ 
 &&{\it orthogonality})&\\ \cline{1-4}
$  C_{4} $ & Coefficients of (1.1), &  Lemma 3.9 &  Positive \\
& interpolation constants &({\it Quasi-reliability})&\\ \cline{1-4}
$ \gamma_0$      &$ \left \{
  \begin{array}{l l}
    0<\gamma_0<1 & ~~ \text{if $T\in \mathcal T_H,$ refined}\\
    1 & \quad \text{otherwise}
  \end{array} \right.$& Corollary  3.5 &  Positive \\ \hline \hline 
  & {\bf II. Contraction property}& & \\ \hline 
 $ \theta_B$     & $0<\theta_B <1$, to be chosen          &  &  (0,1) \\  \cline{1-2} \cline{4-4}
$ \delta_1$    &$<\frac{\theta_B}{4-\theta_B}$     &  Lemma 3.4  &  (0,1) \\  \cline{1-2} \cline{4-4}
$ \rho_B$    &$(1+\delta_1)(1-\frac{\theta_B}{4})$     &  ({\it volume-estimator} &  (0,1) \\  \cline{1-2} \cline{4-4}
$ \varLambda_1$ & $C_1(1+\frac{1}{\delta_1})$     & {\it reduction }) &  (0,1)\\ \hline
$ \theta_A$     &  $0<\theta_A <1$, to be chosen                                &   & (0,1) \\  \cline{1-2} \cline{4-4}
$ \delta_2$     &$\frac{ \theta_A}{4-2 \theta_A}$ & Lemma  3.6  &  (0,1)\\   \cline{1-2} \cline{4-4}
$ \varLambda_2$ & $C_2(1+\frac{1}{\delta_2})$     & ({\it edge-estimator}&  (0,1)\\   \cline{1-2} \cline{4-4}
$ \rho_A$       &$(1+\frac{1}{\delta_2})(1-\frac{\theta_A}{2})$&   {\it reduction }) &  (0,1) \\  \hline
$C_3$            & Coefficient of (1.1) &                         &Positive \\   \cline{1-2} \cline{4-4}
$\alpha_2$       & $0.5+7C_3^2\varLambda_1 H^2$ &  Lemma 3.6 &  (0,1) \\   \cline{1-2} \cline{4-4}
$ \alpha_3$      & $0.5-7C_3^2\varLambda_1 H^2$ &  ({\it Quasi-} &  (0,1) \\   \cline{1-2} \cline{4-4}
$\delta_3$       & $<$1                              & {\it orthogonality}) & (0,1) \\  \cline{1-2} \cline{4-4}
$\alpha_1$       &$\delta_3+4C_3^2\epsilon^2$                             &&   (0,1) \\  \hline
$\alpha_4$       & $\frac{1}{2}\min \{1, \frac{(1-\rho_A) \alpha_3}{({\varLambda_1+ \varLambda_2 }) C_{\rm rel}}\}$&  &  (0,1) \\   \cline{1-2} \cline{4-4}
$\delta_3$        &$<\alpha_4$&Remark 3.2 &(0,1)\\ \cline{1-2} \cline{4-4}
 $\epsilon$       & $(\alpha_4-\delta_3)/(4 C_3^2 )$ &&Positive\\ \cline{1-2} \cline{4-4}
         $h_2$    & $\epsilon ~{\rm and}~ < \frac{1}{\sqrt{14C_3^2\varLambda_1}} $ &&Positive\\ \cline{1-2} \cline{4-4}
$\alpha_1$      &$<\alpha_4$&& (0,1) \\ \hline
$C_5$& $({\varLambda_1+ \varLambda_2 })/{\alpha_3}$&&Positive\\  \cline{1-2} \cline{4-4}
$\alpha$  & $C_5(1-\alpha_1)$&                                                      &  Positive \\  \cline{1-2} \cline{4-4}
$h_3$  & $\min \{ h_1, h_2, \frac{1}{\sqrt{\beta}} \}$&                                                      &  Positive \\  \cline{1-2} \cline{4-4}
$E$       &$2(C_5 \varLambda_3+C_5\alpha_4 C_{rel})$ &  Theorem 4.1  &  Positive \\ \cline{1-2} \cline{4-4}
$\kappa$  & $<\kappa_0:= \frac{({1-\rho_A-C_5\alpha_4C_{rel}})}{E}$&   &  Positive \\  \cline{1-2} \cline{4-4}
$D$&$({1.5+C_5\alpha_4C_{rel}})({\kappa})$ &( {\it Contraction} &  Positive \\  \cline{1-2} \cline{4-4}
 $ \delta_1$   &$\min\{ \frac{E}{2\beta} ,\frac{\theta_B}{4-\theta_B}\}$ &   {\it property} )&  (0,1) \\ \cline{1-2} \cline{4-4}
$\beta$& $2 \max\big{\{}C_5 \varLambda_3 + C_5 \alpha_4 C_{rel}, \frac{C_5\varLambda_3+C_5\alpha_4C_{rel}+D}{1-\rho_B}\big{\}}$&                                 
                                                                                      &  Positive \\  \cline{1-2} \cline{4-4}
$\rho_1$& $\max \Big{\{ } \frac{C_5(1-\alpha_4)}{\alpha}, \rho_A+C_5\alpha_4C_{rel}+E \kappa,$&    &  (0,1)\\ 
&\qquad  \hspace{2cm}$ \frac{((1+\delta_1)\beta-E/2)}{\beta} \Big \}$&&\\  \cline{1-2} \cline{4-4}
$\rho_2$&$\max \Big\{  \frac{C_5(1-\alpha_4)}{\alpha}, 2+C_5\alpha_4C_{rel}-D\kappa,$  & &  (0,1) \\ 
&\qquad  \hspace{2cm}$ \frac{\rho_B\beta+C_5 \varLambda_3+C_5\alpha_4C_{rel}+D }{\beta} \Big \}$&&\\ \cline{1-2} \cline{4-4}
$\rho$& $\max\{\rho_1,\rho_2\}$&  & (0,1) \\  \hline  \hline
  & {\bf III. Quasi-Optimality~}& & \\ \hline 
$h_4$&  $ \frac{1}{2(7+C_3 \varLambda_1+{\beta \varLambda_1})/({C_{\rm eff}^{-1}+\alpha})}$  & Lemma 4.3  &  Positive \\ \cline{1-3} \cline{4-4}
$ \theta_A$     & $\min\{1,\frac{C_{\rm eff}}{9C_4}\}$&   & (0,1) \\ \cline{1-2} \cline{4-4}
$ \theta_B$     & $({\gamma_0^2})/({8(2+C_4C_1h^2)})$& Theorem  4.5 &  (0,1) \\  \cline{1-2} \cline{4-4} 
$\kappa$& $   \min \{ \kappa_0,  \frac{{C_{\rm eff}}/{8}-C_4\theta_A }{4C_4+ {\varLambda_3}/{2}}\big\}$&  ({\it Quasi-} &  Positive \\  \cline{1-2} \cline{4-4}
$\epsilon$& $\min\{ \epsilon, 1, \frac{1}{4C_4} \}$&{optimality}) &  Positive \\ \cline{1-2} \cline{4-4}
$h_5,h_7$& $ \epsilon $& &  Positive \\ \cline{1-2} \cline{4-4}
$\tau_2^2$& $ \frac{(c_{\rm eff}/8)-\kappa ( 4C_4+ {\varLambda_3}/{2})-C_5 \theta_A}{(5/4)C_7} $&  &  Positive \\ \cline{1-2} \cline{4-4}
$h_6$&  $\epsilon $ and $\min\{ h_2, h_4,\frac{1}{\sqrt{4C_4 \varLambda_1}} \}$&  &  Positive \\ \cline{1-2} \cline{4-4}
$\tau_3^2$& $\frac{\beta \gamma_0^2}{8(2+C_4C_1h^2)C_8},C_8=C_4C_1(C_{\rm rel +1})(1+\frac{1}{\kappa}) $ &  &  Positive \\ \cline{1-2} \cline{4-4}
$h_8$&  $\epsilon $ and $\min\{ h_2,h_4, \sqrt{\frac{\alpha}{\beta C_1C_4 }},\sqrt{\frac{\gamma_0^2}{4C_8 }} \}$&  &  Positive \\ \cline{1-2} \cline{4-4}
$h_0$& $\min\{ h_7, h_8 \}$& &  Positive \\ \cline{1-2} \cline{4-4}
$\epsilon_1$& $\tau_1 \xi_\ell $, where $\tau_1:=\min\{ \tau_2, \tau_3 \}$&   &  Positive 
\end{longtable}
\end{center} }

\end{document}